\documentclass{amsart}
\usepackage{amsmath,amscd,amssymb,amsthm,graphicx}

\def\ec{equivalence class}
\def\svc{SVC}

\def\gd{{\gamma, \delta}}
\def\gdp{{\gamma', \delta'}}
\def\tgfd{{\t\gamma_5, \delta}}
\def\lu{{\lambda, \mu}}
\def\lup{{\lambda', \mu'}}
\def\ul{{\mu, \lambda}}
\def\ulp{{\mu', \lambda'}}
\def\ll{{\lambda, \lambda}}
\def\VS{\Vec(S^1)}
\def\EC{\mathop{\rm EC}\nolimits}
\def\LST{\mathop{\rm lower}\nolimits}
\def\UST{\mathop{\rm upper}\nolimits}

\newtheorem{lemma}{Lemma}[section]

\newtheorem{prop}[lemma]{Proposition}
\newtheorem{thm}[lemma]{Theorem}
\newtheorem{cor}[lemma]{Corollary}

\def\cal{\mathcal}

\def\D{{\cal D}} 
 
\def\F{{\cal F}}

\def\frak{\mathfrak}

\def\da{{\frak a}}

\def\db{{\frak b}}

\def\dl{{\frak l}}

\def\ds{{\frak s}}
\def\dsl{\ds\dl}

\def\Bbb{\mathbb}

\def\bC{\Bbb C}

\def\bN{\Bbb N}

\def\bR{\Bbb R}

\def\bZ{\Bbb Z}

\def\ep{\epsilon}

\def\o{\overline}

\def\t{\tilde}

\def\and{\mbox{\rm \ and\ }}

\def\Bol{\mathop{\rm Bol}\nolimits}

\def\Hom{\mathop{\rm Hom}\nolimits}

\def\lac{\mathop{\rm lac}\nolimits}

\def\oh{{\ts\frac{1}{2}}}

\def\PQ{\mathop{\rm PQ}\nolimits}
\def\px{\partial_x}

\def\scrm{\scriptsize\rm}

\def\Span{\mathop{\rm Span}\nolimits}

\def\SQ{\mathop{\rm SQ}\nolimits}

\def\thup{{\mbox{\scrm th}}}

\def\ts{\textstyle}
\def\Vec{\mathop{\rm Vect}\nolimits}

\def\VR{\Vec(\Bbb R)}
\def\VRM{\Vec(\Bbb R^m)}

\def\bgno{\bigbreak\noindent}

\def\cs{composition series}
\def\dog{differential operator}

\def\ie{{\em i.e.,\/}}
\def\iff{if and only if}

\def\irr{irreducible}

\def\meno{\medbreak\noindent}

\def\psdog{pseudo\dog}

\def\r{representation}  

\def\sq{subquotient}
\def\tdm{tensor density module}

\def\th{\thinspace}

\def\vf{vector field}

\title[Equivalence of subquotients]{Equivalence classes of subquotients of pseudodifferential operator modules}

\author{Charles H.\ Conley}
\address{Department of Mathematics \\University of North Texas \\Denton TX 76203, USA} 
\email{conley@unt.edu}

\author{Jeannette M.\ Larsen}
\address{Department of Mathematics \\University of North Texas \\Denton TX 76203, USA} 
\email{JeannetteLarsen@my.unt.edu}

\thanks{The first author was partially supported by Simons Foundation Collaboration Grant 207736.}

\begin{document}

\begin{abstract}
Consider the spaces of \psdog s between \tdm s over the line as modules of the Lie algebra of \vf s on the line.  We compute the \ec es of various \sq s of these modules.  There is a 2-parameter family of \sq s with any given Jordan-H\"older \cs.  In the critical case of \sq s of length~5, the \ec es within each non-resonant 2-parameter family are specified by the intersections of a pencil of conics with a pencil of cubics.  In the cases of resonant \sq s of length~4 with self-dual \cs, as well as of lacunary \sq s of lengths~3 and~4, equivalence is specified by a single pencil of conics.  Non-resonant \sq s of length exceeding~7 admit no non-obvious equivalences.  The cases of lengths~6 and~7 are unresolved.
\end{abstract}


\maketitle

\section{Introduction}  \label{Intro}

Let $\VRM$ be the Lie algebra of polynomial vector fields on $\bR^m$.  Its natural module $\bC[x_1, \ldots, x_m]$ has a 1-parameter family of deformations, the {\em \tdm s\/} $\F_\lambda(\bR^m)$, the sections of the $\lambda^\thup$ power of the determinant bundle.  These modules are defined for all $\lambda \in \bC$, and we write them as
$$ \F_\lambda(\bR^m) := dx^\lambda \bC[x_1, \ldots, x_m]. $$
The Lie action $L_\lambda$ of $\VRM$ on $\F_\lambda(\bR^m)$ is
$$ L_\lambda(X) \bigl(dx^\lambda f(x)\bigr) :=
   dx^\lambda \bigl( X(f) + \lambda f \th \nabla \cdot X \bigr). $$

Let $\D_\lu(\bR^m)$ be the space of \dog s from $\F_\lambda(\bR^m)$ to $\F_\mu(\bR^m)$, and let $\D^k_\lu(\bR^m)$ be the subspace of operators of order~$\le k$.  The Lie action $L_\lu$ of $\VRM$ on $\D_\lu(\bR^m)$ is
$$ L_\lu(X)(T) := L_\mu(X) \circ T - T \circ L_\lambda(X). $$
It preserves the order filtration, and so one has the \sq\ modules
\begin{equation*}
   \SQ^{k,l}_\lu(\bR^m) := \D^k_\lu(\bR^m) / \D^{k-l}_\lu(\bR^m).
\end{equation*}
Note that for $l=1$ these are modules of $k^\thup$ order symbols.  Symbol modules are usually \irr, so we may think of $l$ as the Jordan-H\"older length of $\SQ^{k,l}_\lu(\bR^m)$.  For $l \ge 2$, the \sq\ is usually not completely reducible.

The subject of the present article is the following ``equivalence question'': when are two such \sq s equivalent as modules of $\VRM$?  This question was first posed by Duval and Ovsienko \cite{DO97}, who answered it for modules of the form $\D^k_\ll(\bR^m)$ with $k \le 2$.  In fact they treated smooth \dog s over arbitrary oriented manifolds, but it is a general phenomenon that the result depends only on the dimension of the manifold, and in the Euclidean case it is the same whether one considers polynomial or smooth functions.

Duval and Ovsienko observed a dichotomy between the 1-dimensional and multidimensional cases, due to the fact that for $m \ge 2$, the Jordan-H\"older \cs\ of $\SQ^{k,l}_\lu(\bR^m)$ is determined by the three parameters $\mu-\lambda$, $k$, and $l$, while for $m=1$ it is determined only by the two parameters $\mu-\lambda-k$ and $l$.  This in turn is because the symbol modules of $\D_{\lu}(\bR^m)$ are \tdm s for $m=1$, but not for $m \ge 2$.

The work \cite{DO97} inspired several articles.  (We note that some authors use a different sign convention, writing $\F_{-\lambda}$ where we write $\F_\lambda$.)  In the multidimensional case $m \ge 2$, Lecomte, Mathonet, and Tousset \cite{LMT96} determined the \ec es of the modules $\D^k_\ll(\bR^m)$ for $k \ge 3$, and Gargoubi and Ovsienko \cite{GO96} did the same in the 1-dimensional case.  Genuine \sq s were first considered by Lecomte and Ovsienko \cite{LO99}, who also made the natural and important generalization to \psdog s in the 1-dimensional case, allowing the order $k$ to vary continuously.  They computed the \ec es of the modules $\SQ^{k,l}_\ll(\bR^m)$ (but only for generic values of $k$ when $m=1$).

The equivalence question for arbitrary $(\lu)$ was first considered by Mathonet \cite{Ma99} and Gargoubi \cite{Ga00}.  In \cite{Ma99}, all $\VRM$-intertwining maps between the modules $\D^k_\lu(\bR^m)$ are determined for $m \ge 2$.  To our knowledge, the question is not yet settled in the multidimensional case for genuine \sq s $\SQ^{k,l}_\lu(\bR^m)$.  The classification of the $\VRM$-maps from \dog\ modules $\D^k_\lu(\bR^m)$ to \tdm s $\F_\nu(\bR^m)$, carried out in \cite{Ma00}, is closely related.

The equivalences among the modules $\D^k_\lu(\bR)$ are determined in \cite{Ga00}.  This article contains errors pointed out in \cite{CS04} which call the results at $k=3$ and $\mu-\lambda=1$ or~$3$ into question, but in fact as we shall see here they are correct.

In all of these articles on the equivalence question there is a 1-parameter family of modules with any given composition series.  The various answers obtained have a common trait: for small lengths~$l$, modules with the same composition series are all equivalent excepting a finite number of special cases, while for larger lengths~$l$, modules are equivalent only to their conjugates (adjoints).  The most interesting cases involve the critical intermediate lengths, which in these articles are always $l=3$ or $l=4$.

In this article we consider the \ec es of the modules $\SQ^{k,l}_\lu(\bR)$, the most general 1-dimensional setting.  Here there is a two parameter family of modules with any given \cs, which causes the critical length to increase to $l=5$.  At this length there is a new phenomenon: for each \cs, the \ec es are generically six pairs of conjugate modules, determined by the intersections of a certain pencil of conics with a certain pencil of cubics.  During the course of the analysis we obtain results in lengths $l \le 4$ unifying the 1-dimensional results of \cite{DO97}, \cite{GO96}, \cite{LO99}, \cite{Ga00}, and \cite{Ma00}.

We have been unable to resolve the case of length~6.  We can reduce the equivalence question to the computation of a certain Gr\"obner basis, but the standard software packages were unable to find this basis on the computers available to us.  We expect that there are only the obvious equivalences arising from conjugation and the de Rham differential; exceptional equivalences would be interesting.  We have not resolved the case of length~7 either, but this will be much easier because one has six invariants rather than four.  We can prove that in lengths $l \ge 8$ there are no non-obvious equivalences, but we have not included the details here.

We also study ``lacunary \sq s'', $\VR$-modules whose \cs\ are missing certain symbol modules.  For example, the lacunary \sq s of length~3 composed of the order~$k$, $k-2$, and $k-4$ symbols have the same two parameters for each \cs, and generically they are equivalent \iff\ their parameters lie on the same member of the pencil of conics involved in the $l=5$ case.  On the other hand, the \ec es of the lacunary \sq s of length~4 composed of the order~$k$, $k-2$, $k-3$, and $k-5$ symbols are determined by a new pencil of conics.

As in \cite{GO96}, \cite{LO99}, and \cite{Ga00}, our main tool is the {\em projective quantization,\/} the decomposition of $\D_\lu(\bR)$ under the action of the {\em projective subalgebra\/} of $\VR$, a copy of $\dsl_2$.  Generically one has complete reducibility under this subalgebra, in which case it suffices to use the formulas for the action of $\VR$ with respect to the projective quantization deduced by Cohen-Tretkoff, Manin, and Zagier \cite{CMZ97}.  The exceptions are the {\em resonant\/} cases, where we use the modifications of these formulas obtained in \cite{Ga00} and \cite{CS04}.  In fact we review these resonant formulas in more detail than is needed for the equivalence question in order to give them in a much simpler form; see Theorem~\ref{res CMZ}.  One consequence of this simplification is Corollary~\ref{B factors}, which explains certain initially mysterious factorizations of the non-resonant formulas.

The length~4 resonant case $\SQ^{\mu-\lambda+1,4}_\lu(\bR)$ with self-dual \cs\ is particularly interesting.  The special case $\D^3_{\lambda,\lambda+2}(\bR)$ of \dog s was studied in \cite{Ga00}: it is the only length $l \le 4$ \cs\ for which there are no equivalences other than conjugation.  Here we see that in the more general setting of \sq s, equivalence is determined by a single pencil of conics, not one of the pencils of conics arising at $l=5$.

A preliminary outline of these results was given in \cite{Co09}, and the non-resonant cases comprise the Ph.D.\ thesis \cite{La12} of the second author.  The content of the article is as follows.  In Section~\ref{Dfns} we state the equivalence question for \psdog s and recall conjugation, the Adler trace, the de Rham differential, and resonance.  In Section~\ref{Nonres Results} we state the complete answer to the equivalence question in all non-resonant cases of length $l \le 5$, and in Section~\ref{Res Results} we do the same in all resonant cases of length $l \le 4$.  In Section~\ref{Lac Results} we state our results on the lacunary equivalence question, and in Section~\ref{Pictures} we discuss the various pencils of conics and cubics which arise.  All proofs are given in Section~\ref{Proofs}, and we conclude in Section~\ref{Remarks} with remarks on the equivalence question in higher lengths.  Let us mention Proposition~\ref{noncocycles}, which states that non-resonant \sq s of length $l \ge 6$ are equivalent \iff\ each of their own length~5 \sq s are equivalent, except possibly in certain cases involving the ``$\sqrt{19}$'' 1-cocycles of $\VR$ discovered by Feigin and Fuchs \cite{FF80}.

\section{Definitions and background}  \label{Dfns} 

Henceforth we work exclusively in one dimension, so we will drop the argument $\bR$ and write simply $\F_\lambda$ for $dx^\lambda \bC[x]$, $\D_\lu$ for the \dog s from $\F_\lambda$ to $\F_\mu$, and so on.  We adopt the standard convention of writing $\delta$ for $\mu-\lambda$:
\begin{equation*}
   \delta(\lambda,\mu) := \mu - \lambda.
\end{equation*}
We denote the non-negative integers by $\bN$ and the positive integers by $\bZ^+$.

For any $k \in \bC$, the $\VR$-module of \psdog s from $\F_\lambda$ to $\F_\mu$ of order in $k-\bN$ consists of formal sums:
\begin{equation*}
   \Psi^k_\lu := \Bigl\{ dx^\delta \sum_{i=0}^\infty f_i(x) \px^{k-i}: f_i \in \bC[x] \Bigr\}.
\end{equation*}
The action $L_\lu$ of $\VR$ on $\Psi^k_\lu$ is the natural extension of the action on $\D_\lu$:
\begin{equation*}
   L_\lu(g \px) (dx^\delta f \px^r) := 
   dx^\delta \Bigl\{\bigl( gf' + (\delta - r)g'f \bigr) \px^r -
   f \sum_{s=1}^\infty {r \choose s} \Bigl(\lambda + \frac{r-s}{s+1} \Bigr)
   g^{(s+1)} \px^{r-s} \Bigr\}.   
\end{equation*}

Observe that $\Psi^{k-1}_\lu$ is a submodule of $\Psi^k_\lu$, and $dx^\delta f \px^k \mapsto dx^{\delta-k} f$ defines a $\VR$-equivalence from $\Psi^k_\lu/\Psi^{k-1}_\lu$ to $\F_{\delta-k}$.  We extend the definition of $\SQ^{k,l}_\lu$ to \psdog s: for $\lambda$, $\mu$, and~$k$ in $\bC$ and $l \in \bZ^+$, set
\begin{equation*}
   \SQ^{k,l}_\lu := \Psi^k_\lu / \Psi^{k-l}_\lu.
\end{equation*}
We shall refer to this $\VR$-module as a \sq\ of length~$l$ with \cs\ $\{ \F_{\delta-k}, \F_{\delta-k+1}, \ldots, \F_{\delta-k+l-1} \}$.  This is slightly inaccurate: although $\F_\nu$ is \irr\ for $\nu \not= 0$, $\F_0$ is a module of length~2 with \cs\ $\{\F_1, \bC\}$.  However, the following lemma is clear.

\begin{lemma} \label{n}
If\/ $\SQ^{k,l}_\lu$ and\/ $\SQ^{k',l'}_\lup$ are equivalent, then $l=l'$ and $\delta-k = \delta'-k'$.
\end{lemma}

We may state the equivalence question as follows.  Define
$$ n(k,\delta) := \delta-k. $$

\meno{\bf Question.}
For fixed $n\in\bC$ and $l\in\bZ^+$, what are the $\VR$-\ec es of the set
\begin{equation*}
   \bigl\{ \SQ^{\delta-n,l}_\lu: \lu \in \bC \bigr\}
\end{equation*}
of \sq s of length~$l$ with \cs\ $\{ \F_n, \F_{n+1}, \ldots, \F_{n+l-1} \}$?

\medbreak
We should note that this question contains the equivalence question for the \dog\ modules themselves.  Indeed, for $k\in\bN$ the length $k+1$ \sq\ $\SQ^{k,k+1}_\lu$ is simply $\D^k_\lu$.  More generally, for $l \ge k+2$ we have the canonical $\VR$-splittings
\begin{equation} \label{Dsplitting}
   \Psi^k_\lu = \D^k_\lu \oplus \Psi^{-1}_\lu, \qquad
   \SQ^{k,l}_\lu = \D^k_\lu \oplus \SQ^{-1,l-k-1}_\lu.
\end{equation}

Our main result is the answer to the equivalence question in all cases with $l \le 4$, and in the non-resonant cases with $l = 5$.  In order to state it efficiently, we recall conjugation, the Adler trace, the de Rham differential, and resonance.

{\em Conjugation of \psdog s\/} is the adjoint map $T \mapsto T^*$ from $\Psi^k_\lu$ to $\Psi^k_{1-\mu,1-\lambda}$ defined by
\begin{equation*}
   (dx^\delta f \px^r)^* := e^{\pi ir} dx^\delta \px^r f
   = e^{\pi ir} \sum_{s=0}^\infty {r \choose s} f^{(s)} \px^{r-s}.
\end{equation*}
Observe that conjugating twice acts on $\Psi^k_\lu$ as the scalar map $e^{2\pi i k}$.  The following lemma is well-known and easy to prove.

\begin{lemma} \label{conj}
Conjugation is a $\VR$-equivalence from $\Psi^k_\lu$ to $\Psi^k_{1-\mu,1-\lambda}$.  In particular, as $\VR$-modules,
$$ \SQ^{\delta-n,l}_\lu \cong \SQ^{\delta-n,l}_{1-\mu,1-\lambda}. $$
\end{lemma}

The {\em Adler trace,\/} also known as the {\em noncommutative residue,\/} exists in the category of $\VS$-modules.  Algebraically, one passes to this category by simply adjoining $x^{-1}$ to $\VR$ and all of its modules.  The trace is a nondegenerate $\VS$-invariant pairing between $\Psi^{\bN+k}_\lu(S^1)$ and $\Psi^{\bN-k}_\ul(S^1)$; see \cite{CMZ97, LO99, CS04}.  It yields the following lemma.

\begin{lemma} \label{trace}
The $\VS$-modules\/ $\SQ^{k,l}_\lu(S^1)$ and\/ $\SQ^{l-2-k,l}_\ul(S^1)$ are dual.
\end{lemma}

Lemmas~\ref{conj} and~\ref{trace} both originate in the fact that the $\VS$-modules $\F_\nu(S^1)$ and $\F_{1-\nu}(S^1)$ are dual.  The consequence of Lemma~\ref{trace} relevant here is the following corollary, given as Lemma~4.2 in \cite{Co05}.

\begin{cor} \label{dual}
As\/ $\VR$-modules,
$$ \SQ^{\mu-\lambda-n,l}_{\lu} \cong \SQ^{\mu'-\lambda'-n,l}_\lup
   \mbox{\it\ \ \iff\ \ }
   \SQ^{l-2+n+\lambda-\mu,l}_\ul \cong \SQ^{l-2+n+\lambda'-\mu',l}_\ulp. $$
\end{cor}

The following definition will permit us to take advantage of these symmetries.

\meno {\bf Definition.}
Let $\gamma(\lambda,\mu) := 3(\lambda+\mu-1)^2$ and $N_l(n) := n + \oh l -1$.

\medbreak
For fixed $l$ and $n$, Lemma~\ref{conj} implies that the \ec\ of $\SQ^{\delta-n,l}_\lu$ depends only on $\gamma$ and $\delta$.  By Corollary~\ref{dual}, the equations defining this \ec\ are symmetric under $(N_l,\gd) \mapsto (-N_l,\gamma,-\delta)$.  Therefore we will give the equations in terms of these coordinates.  

Keep in mind that $(\gd)$ specifies a conjugate pair of values of $(\lu)$ rather than a single value.  In fact, many of our formulas involve $\gamma^{1/2}$.  Although the statements of our main results are independent of the choice of sign of the square root, for concreteness we specify
\begin{equation*}
   \gamma^{1/2}(\lu) := \sqrt{3}\, (\lambda+\mu-1).
\end{equation*}

Henceforth we will always use the notation
\begin{equation*}
   (\gdp) := \bigl( \gamma(\lup), \delta(\lup) \bigr).
\end{equation*}
We make the following definition in order to be able to regard the \ec\ of $\SQ^{\delta-n,l}_\lu$ as a subset of the $(\gd)$-plane.

\meno {\bf Definition.}
$\EC^l_n(\gamma,\delta) := \bigl\{ (\gamma',\delta') \in \bC^2:\ 
\SQ^{\delta'-n,l}_\lup \cong \SQ^{\delta-n,l}_\lu \bigr\}$.

\medbreak
The {\em de Rham differential\/} is
$$ d := dx \th \px:\F_0 \to \F_1, $$
the only non-scalar $\VR$-map between \tdm s.  It gives rise to an equivalence between \sq s of arbitrary length which we now describe.  Write $L_d$ and $R_d$ for left and right composition with $d$, respectively:
\begin{equation*} \begin{array}{rl}
   L_d: \Psi^k_{\lambda,0} \to \Psi^{k+1}_{\lambda,1}, &
   T \mapsto d \circ T, \\[6pt]
   R_d: \Psi^k_{1,\mu} \to \Psi^{k+1}_{0,\mu}, &
   T \mapsto T \circ d.
\end{array} \end{equation*}
These maps are both $\VR$-isomorphisms, which induce $\VR$-isomorphisms
\begin{equation*}
   L_d: \SQ^{-\lambda-n,l}_{\lambda,0} \to \SQ^{1-\lambda-n,l}_{\lambda,1}, \qquad
   R_d: \SQ^{\mu-1-n,l}_{1,\mu} \to \SQ^{\mu-n,l}_{0,\mu}.
\end{equation*}
Observe that the two cases form a conjugate pair, so in $(\gd)$-coordinates they appear as a single case.  Thus we have:

\begin{lemma} \label{Bol}
For all $l$, $n$, and $\nu$, $\EC^l_n \bigl( 3(\nu+1)^2, \nu \bigr) = \EC^l_n \bigl( 3\nu^2, \nu+1 \bigr)$.
\end{lemma}

The interplay between~(\ref{Dsplitting}) and the maps $L_d$ and $R_d$ gives:

\begin{lemma} \label{DBol}
For all $\lambda$ and $\mu$, we have the\/ $\VR$-splittings
\begin{equation*} \begin{array}{rclrcl}
   \D_{\lambda,1} &=& \D^0_{\lambda,1} \oplus L_d(\D_{\lambda,0}), \quad & \quad
   \Psi^{-1}_{\lambda,0} &=& \Psi^{-2}_{\lambda,0} \oplus L_d^{-1}(\D^0_{\lambda,1}), \\[6pt]
   \D_{0,\mu} &=& \D^0_{0,\mu} \oplus R_d(\D_{1,\mu}), \quad & \quad
   \Psi^{-1}_{1,\mu} &=& \Psi^{-2}_{1,\mu} \oplus R_d^{-1}(\D^0_{0,\mu}), \\[6pt]
   \D_{0,1} &=& \D^1_{0,1} \oplus L_d R_d (\D_{1,0}), &
   \Psi^{-1}_{1,0} &=& \Psi^{-3}_{1,0} \oplus L_d^{-1} R_d^{-1} (\D^1_{0,1}).
\end{array} \end{equation*}
\end{lemma}

{\em Resonance\/} is the failure of complete reducibility under the action of the {\em projective subalgebra\/} of $\VR$.  This subalgebra is
\begin{equation*}
   \da := \Span_\bC \bigl\{ \px, x\px, x^2\px \bigr\},
\end{equation*}
an isomorphic copy of $\dsl_2$.  Its Casimir operator,
\begin{equation*}
   Q := (x\px)^2 - (x\px) - (x^2\px)(\px),
\end{equation*}
acts on $\F_\nu$ by the scalar $L_\nu(Q) = \nu^2-\nu$.

The module $\SQ^{\delta-n,l}_\lu$ cannot be resonant unless its composition series has repeated Casimir eigenvalues.  Since the values of $\nu^2-\nu$ are symmetric around $\nu = \oh$, this occurs \iff\ $(n+i) + (n+j)$ is $1$ for some $0 \le i < j \le l-1$, \ie\ $l \ge 2$ and $-2n \in \{0, 1, \ldots, 2l-4\}$.  In fact, $\SQ^{\delta-n,l}_\lu$ is generically resonant for such values of $n$, so we make the following definition.

\meno{\bf Definition.}
$\SQ^{\delta-n,l}_\lu$ is {\em resonant\/} if $n \in \bigl\{0, -\oh, -1, -\frac{3}{2}, \ldots, 2 - l \bigr\}$.  We say that it is {\em integral resonant\/} or {\em half-integral resonant\/} depending on whether $n$ is integral or half-integral.

\medbreak
Note that the set of resonant values of $N_l$ is $\bigl\{0, \pm\oh, \pm 1, \pm\frac{3}{2}, \ldots, \pm(\oh l -1) \bigr\}$.  Its symmetry around zero is a consequence of Lemma~\ref{trace}.  Resonant \psdog\ modules were studied in detail in \cite{Ga00} and \cite{CS04}.

\section{Non-resonant results}  \label{Nonres Results}

In this section we answer the equivalence question in all non-resonant cases of length $l \le 5$.  Proofs will be deferred to Section~\ref{Proofs}.  Keep in mind that by Lemma~\ref{n}, $l$ and $n$ are invariant under equivalence: they are complete invariants for the \cs\ of $\SQ^{\delta-n,l}_\lu$.

For $l=1$ there is nothing to prove: there are no resonant cases and $\SQ^{\delta-n,1}_\lu$ is always equivalent to $\F_n$.  For $l=2$, the only resonant case is $n=0$.  For $n \not= 0$ it is well-known that there is only one \ec: for all $(\lu)$,
$$ \SQ^{\delta-n,2}_\lu \cong \F_n \oplus \F_{n+1}. $$
See for example Lemma~7.9 of \cite{LO99}, Lemma~3.3 of \cite{CS04}, or page~72 of \cite{Co09}.

In order to state the results for $3 \le l \le 5$, we give slight modifications of~(6), (7), and~(8) of \cite{Co05} and make a convenient definition:
\begin{equation} \label{3Bs} \begin{array}{rcl}
   B_{m+2,m}(\gd) &:=& \gamma - \bigl[ (2m+1) \delta + (m^2+m+1) \bigr], \\[6pt]
   B_{m+3,m}(\gd) &:=& \gamma^{3/2} - 3\gamma^{1/2} \bigl[(m+1) \delta + 1 \bigr], \\[6pt]
   B_{m+4,m}(\gd) &:=& \gamma^2 - \gamma \bigl[2(2m+3) \delta - (2m^2+6m-3)\bigr] \\[6pt]
               && - {\ts\frac{3}{5}} m(m+3) \bigl[2(2m+3) \delta + (m^2+3m+6) \bigr].
\end{array} \end{equation}

\bgno {\bf Definition.}
Two \sq s $\SQ^{\delta-n,l}_\lu$ and $\SQ^{\delta'-n,l}_\lup$ are said to {\em induce simultaneous vanishing of the functions\/} $f_1(\gd), \ldots, f_r(\gd)$ if for all~$s$, $f_s(\gd)$ and $f_s(\gdp)$ are either both zero or both non-zero.

\medbreak
Henceforth we will use the Pochhammer symbol $(x)_r$ for the falling factorial:
$$ (x)_r := x(x-1)\cdots (x-r+1), \qquad (x)_0 := 1. $$

\subsection{Length $l=3$} \label{l3}

Here the set of resonant values of $n$ is $\bigl\{-1, -\oh, 0 \bigr\}$, and so that of $N_3 = n+\oh$ is $\bigl\{0, \pm\oh \bigr\}$.  We will need the formula for $B_{n+2,n}$ in terms of $N_3$:
$$ B_{n+2,n} = \gamma - \bigl[ 2N_3 \delta + N_3^2 + {\ts\frac{3}{4}} \bigr]. $$

\begin{prop} \label{nonres3}
For $n$ non-resonant, $\SQ^{\delta-n,3}_\lu$ and $\SQ^{\delta'-n,3}_\lup$ are equivalent \iff\ they induce simultaneous vanishing of 
\begin{equation} \label{l3poly}
   (\delta-N_3+\oh)_2\, B_{n+2,n}.
\end{equation}
The \ec\ where~(\ref{l3poly}) vanishes splits as $\F_n \oplus \F_{n+1} \oplus \F_{n+2}$.
\end{prop}

Let us make some remarks on this proposition and explain how to use it to recover earlier length~3 results.  First, $(\delta-N_3\pm\oh)_2 = 0$ \iff\ the order $\delta-n$ is~$0$ or~$1$.  In these cases the \ec\ is split because of~(\ref{Dsplitting}).

Observe that~(\ref{l3poly}) has the symmetry promised by Corollary~\ref{dual}.  To verify Lemma~\ref{Bol} directly, check that the two $(\gd)$-values $\bigl( 3(\nu+1)^2, \nu \bigr)$ and $\bigl( 3\nu^2, \nu+1 \bigr)$ give the same value in~(\ref{l3poly}), namely,
\begin{equation*}
   (\nu-N_3+{\ts\frac{3}{2}})_3\, (3\nu + N_3 + {\ts\frac{3}{2}}).
\end{equation*}

Before recapitulating earlier results we define
\begin{equation*}
   H_k(\lu) := 4 (3\lambda+k-2) (3\mu-k-1) + (k-2)(k+1).
\end{equation*}
The zeroes of this function form the hyperbola given in Definition~3.3 of \cite{Ma00}.  The reader may check that $H_k(\lu) = 3B_{n+2,n}(\gd)$.

\cite{DO97} gives the \ec es of the modules $\D^2_\ll = \SQ^{2,3}_\ll$: the $\lambda$-values~$0$ and~$1$ form the split class and all the other $\lambda$-values form the other class.  To recover this result, note that at $n = -2$, $\delta=0$, and $\gamma=3(2\lambda-1)^2$ we have
$$ B_{0,-2}(\gamma,0) = {\ts\frac{1}{3}} H_2(\ll) = 12\lambda(\lambda-1). $$

\cite{LO99} gives the \ec es of the modules $\SQ^{k,3}_\ll$ with $k\not\in\oh\bN$: the two $\lambda$-roots of $12\lambda(\lambda-1) - (k-2)(k+1)$ form the split class and all other $\lambda$-values form the other class.  This is because at $n=-k$ we have
\begin{equation*}
   B_{2-k,-k}(\gamma,0) = {\ts\frac{1}{3}} H_k(\ll) =
   12\lambda(\lambda-1) - (k-2)(k+1).
\end{equation*}

\cite{Ga00} gives the \ec es of the modules $\D^2_\lu = \SQ^{2,3}_\lu$: modules with a given $\delta$ have two classes, the conjugate pair $\bigl\{ \D^2_{0,\delta}, \D^2_{1-\delta,1} \bigr\}$ and all the others.  This result is a manifestation of Lemma~\ref{DBol}.  To recover it in the non-resonant cases, note that at $k=2$ we have $n=\delta-2$ and
$$ B_{\delta,\delta-2}(\gd) = {\ts\frac{1}{3}} H_2(\lu) = 12\lambda(\mu-1). $$

\cite{Ma00} gives all examples of projections $\D^k_\lu \twoheadrightarrow \F_{\mu-\lambda-k+2}$: they exist \iff\ $k \in 2+\bN$ and $H_k(\lu) = 0$.  The explanation is that in the non-resonant case, such projections can only exist when $\SQ^{k,3}_\lu$ is in the split \ec.

\subsection{Length $l=4$} \label{l4}

Here the set of resonant values of $n$ is $\{-2, -\frac{3}{2}, -1, -\oh, 0\}$, and so that of $N_4 = n+1$ is $\{0, \pm\oh, \pm 1\}$.  In terms of $N_4$,
\begin{eqnarray*}
   B_{n+2,n} &=& \gamma - \bigl[ (2N_4-1) \delta
   + (N_4^2-N_4+1) \bigr], \\[6pt]
   B_{n+3,n+1} &=& \gamma - \bigl[ (2N_4+1) \delta
   + (N_4^2+N_4+1) \bigr], \\[6pt]
   B_{n+3,n} &=& \gamma^{3/2} - 3\gamma^{1/2} \bigl[ N_4 \delta + 1 \bigr].
\end{eqnarray*}

\begin{prop} \label{nonres4}
For $n$ non-resonant, $\SQ^{\delta-n,4}_\lu$ and $\SQ^{\delta'-n,4}_\lup$ are equivalent \iff\ they induce simultaneous vanishing of 
\begin{equation*}
   (\delta-N_4+1)_2\, B_{n+2,n}, \quad
   (\delta-N_4)_2\, B_{n+3,n+1}, \quad
   (\delta-N_4+1)_3\, B_{n+3,n}.
\end{equation*}
\end{prop}

Note that all eight possible sets of vanishings occur, so there are eight \ec es.  However, all three polynomials vanish only when there is vanishing among the Pochhammer symbols.  This occurs when the order $\delta-n$ is $0$, $1$, or~$2$, the situation of~(\ref{Dsplitting}).  For example, all modules with $\delta = N_4$ are in this \ec.  We now use Proposition~\ref{nonres4} to recover earlier length~4 results.  

\cite{GO96} computes the \ec es of $\D^3_\ll = \SQ^{3,4}_\ll$: there are four, given by the sets of $\lambda$-values
$$ \bigl\{\oh\bigr\}, \quad \bigl\{0,1\bigr\}, \quad
   \bigl\{\oh \Bigl(1\pm\sqrt{\ts\frac{7}{3}}\th \Bigr) \bigr\}, \quad
   \bC \setminus \bigl\{ 0, \oh, 1, \oh \Bigl(1\pm\sqrt{\ts\frac{7}{3}}\th \Bigr) \bigr\}. $$  To prove this, apply the proposition to the following equalities:
\begin{eqnarray*}
   B_{-1,-3}(\gamma,0) &=& {\ts\frac{1}{3}} H_3(\ll) \ =\ 4(3\lambda^2-3\lambda-1), \\[6pt]
   B_{-2,0}(\gamma,0) &=& {\ts\frac{1}{3}} H_2(\ll) \ =\ 12\lambda(\lambda-1), \\[6pt]
   B_{0,-3}(\gamma,0) &=& 24\sqrt{3}\th \lambda (\lambda-\oh) (\lambda-1).
\end{eqnarray*}

\cite{LO99} gives the generic \ec\ of the modules $\SQ^{k,4}_\ll$ with $k\not\in\oh\bN$: all those on which none of the functions $t^2_k$, $t^2_{k-1}$, and $t^3_k$ defined in Proposition~7.10 of that article vanish are equivalent.  To prove this, check that these three functions are proportional to
\begin{eqnarray*}
   B_{2-k,-k}(\gamma,0) &=& {\ts\frac{1}{3}}
   H_k(\ll) \ =\ 12\lambda (\lambda-1) - (k-2)(k+1), \\[6pt]
   B_{3-k,1-k}(\gamma,0) &=& {\ts\frac{1}{3}}
   H_{k-1}(\ll) \ =\ 12\lambda (\lambda-1) - k(k-3), \\[6pt]
   B_{3-k,-k}(\gamma,0) &=& 24\sqrt{3}\th \lambda (\lambda-\oh) (\lambda-1).
\end{eqnarray*}

\cite{Ga00} gives the \ec es of the modules $\D^3_\lu = \SQ^{3,4}_\lu$, superseding \cite{GO96}.  To reconstruct the results in the non-resonant cases, note that here $n=\delta-3$ and the Pochhammer symbols in Proposition~\ref{nonres4} never vanish, so the \ec es are determined by the vanishing of
\begin{eqnarray*}
   B_{\delta-1,\delta-3}(\gd) &=&
   {\ts\frac{1}{3}} H_3(\lu) \ =\ {\ts\frac{4}{3}} 
   \bigl( (3\lambda+1) (3\mu-4) + 1 \bigr), \\[6pt]
   B_{\delta,\delta-2}(\gd) &=&
   {\ts\frac{1}{3}} H_2(\lu) \ =\ 12\lambda(\mu-1), \\[6pt]
   B_{\delta,\delta-3}(\gd) &=& 12\sqrt{3}\th \lambda (\mu-1) (\lambda+\mu-1).
\end{eqnarray*}

\cite{Ma00} gives all examples of projections $\D^k_\lu \twoheadrightarrow \F_{\mu-\lambda-k+3}$: they can exist only if $k\in 3+\bN$.  For $k=3$ they exist \iff\ either $\lambda=0$ or $\mu=1$, while for $k \ge 4$ they exist only in the two self-conjugate cases
$$ \lambda = -{\ts\frac{1}{3}} (k-3) \pm {\ts\frac{1}{6}} \sqrt{k(k-3)} \th ,
   \qquad \mu = 1-\lambda. $$
The explanation is that in the non-resonant case, such projections can exist only if both $B_{n+3,n+1} \propto H_{k-1}$ and $B_{n+3,n} \propto \lambda (\mu-1) (\lambda+\mu-1)$ are zero.

\subsection{Length $l=5$} \label{l5}

We have seen that in length $l \le 4$, almost all \sq s $\SQ^{\delta-n,l}_\lu$ with a given~$n$ are equivalent.  At $l=5$ there is a new phenomenon: equivalence is determined by two rational invariants.  We begin with a special case of the definition of simultaneous vanishing:

\meno {\bf Definition.}
Two non-resonant \sq s $\SQ^{\delta-n,l}_\lu$ and $\SQ^{\delta'-n,l}_\lup$ are said to satisfy the {\em simultaneous vanishing condition\/} (\svc) if for all $(i,j)$ such that $0 \le j \le i-2$ and $2 \le i \le l-1$, they induce simultaneous vanishing of
\begin{equation} \label{Bij}
   (\delta-n-j)_{i-j}\, B_{n+i,n+j}(\gd).
\end{equation}

\medbreak
As before, the Pochhammer symbols only vanish in the situation of~(\ref{Dsplitting}).  With this definition we can restate Propositions~\ref{nonres3} and~\ref{nonres4} concisely:

\begin{prop} \label{l2l3l4}
For $l \le 4$, non-resonant \sq s\/ $\SQ^{\delta-n,l}_\lu$ and\/ $\SQ^{\delta'-n,l}_\lup$ are equivalent \iff\ they satisfy the \svc.
\end{prop}

This is not true for $l \ge 5$.  We expect that in length $l \ge 6$, the conjugation and de Rham equivalences described in Lemmas~\ref{conj} and~\ref{Bol} are the only equivalences.  However, in length~5 there are more.  Here the set of resonant values of $n$ is $\{-3, -\frac{5}{2}, -2, -\frac{3}{2}, -1, -\oh, 0\}$, and so that of $N_5 = n+\frac{3}{2}$ is $\{0, \pm\oh, \pm 1, \pm\frac{3}{2} \}$.  In terms of $N_5$,
\begin{equation} \label{6Bijs} \begin{array}{rcl}
   B_{n+2,n} &=& \gamma - \bigl[ 2(N_5-1) \delta
   + (N_5-1)^2 + {\ts\frac{3}{4}} \bigr], \\[6pt]
   B_{n+3,n+1} &=& \gamma - \bigl[ 2N_5 \delta +
   N_5^2 + {\ts\frac{3}{4}} \bigr], \\[6pt]
   B_{n+4,n+2} &=& \gamma - \bigl[ 2(N_5+1) \delta + 
   (N_5+1)^2 + {\ts\frac{3}{4}} \bigr], \\[6pt]
   B_{n+3,n} &=& \gamma^{3/2}
   - 3\gamma^{1/2} \bigl[ (N_5 - \oh) \delta + 1 \bigr], \\[6pt]
   B_{n+4,n+1} &=& \gamma^{3/2} 
   - 3\gamma^{1/2}\bigl[ (N_5 + \oh) \delta + 1 \bigr], \\[6pt]
   B_{n+4,n} &=& \gamma^2 - 2\gamma
   \bigl[ 2 N_5 \delta - N_5^2 + {\ts\frac{15}{4}} \bigr]
   - {\ts\frac{3}{5}} (N_5^2 - {\ts\frac{9}{4}})
   \bigl[ 4 N_5 \delta + N_5^2 + {\ts\frac{15}{4}} \bigr].
\end{array} \end{equation}
The invariants which will determine equivalence are
\begin{eqnarray*}
I_n(\gd) &:=& B_{n+4,n} \big/ B_{n+4,n+2} B_{n+2,n}, \\[6pt]
J_n(\gd) &:=& B_{n+4,n} B_{n+3,n+1} \big/ B_{n+4,n+1} B_{n+3,n}, \\[6pt]
K_n(\gd) &:=& B_{n+4,n+2} B_{n+3,n+1} B_{n+2,n} \big/ B_{n+4,n+1} B_{n+3,n}.
\end{eqnarray*}
Note that $K_n = J_n / I_n$.  We now state our main result in the non-resonant case.

\begin{thm} \label{nonres5}
For $n$ non-resonant, the \sq s\/ $\SQ^{\delta-n,5}_\lu$ and\/ $\SQ^{\delta'-n,5}_\lup$ are equivalent \iff\ they satisfy the \svc, and either $(\delta-n)_4\, (\delta'-n)_4 = 0$, or $(\delta-n)_4\, (\delta'-n)_4 \not= 0$ and one of the following mutually exclusive conditions holds:
\begin{enumerate}
\item[(i)]
At least two of\/ $B_{n+4,n}$, $B_{n+4,n+2}\, B_{n+2,n}$, and $B_{n+4,n+1}\, B_{n+3,n}\, B_{n+3,n+1}$ are zero.
\smallbreak \item[(ii)]
$B_{n+4,n}\, B_{n+4,n+2}\, B_{n+2,n}$ is not zero, $B_{n+4,n+1}\, B_{n+3,n}\, B_{n+3,n+1}$ is zero, and
\begin{equation*}
   I_n(\gd)\ =\ I_n(\gdp).
\end{equation*}
\smallbreak \item[(iii)]
$B_{n+4,n}\, B_{n+4,n+1}\, B_{n+3,n}\, B_{n+3,n+1}$ is not zero, $B_{n+4,n+2}\, B_{n+2,n}$ is zero, and
\begin{equation*}
   J_n(\gd)\ =\ J_n(\gdp).
\end{equation*}
\smallbreak \item[(iv)]
$B_{n+4,n+1}\, B_{n+3,n}\, B_{n+4,n+2}\, B_{n+3,n+1}\, B_{n+2,n}$ is not zero, $B_{n+4,n}$ is zero, and
\begin{equation*}
   K_n(\gd)\ =\ K_n(\gdp).
\end{equation*}
\smallbreak \item[(v)]
$B_{n+4,n}\, B_{n+4,n+1}\, B_{n+3,n}\, B_{n+4,n+2}\, B_{n+3,n+1}\, B_{n+2,n}$ is not zero, and
\begin{equation*}
   I_n(\gd)\ =\ I_n(\gdp), \qquad J_n(\gd)\ =\ J_n(\gdp).
\end{equation*}
\end{enumerate}
\end{thm}

Let us give a preliminary interpretation of this theorem.  Recall that $n$ and thus also $N_5$ may be regarded as fixed because they determine the \cs\ of $\SQ^{\delta-n,5}_\lu$ and so are invariant under equivalence.  By Theorem~\ref{nonres5}(v), away from the zero loci of the six $B_{n+i,n+j}$'s (\ref{6Bijs}) in the $(\gd)$-plane, $I_n$ and $J_n$ are complete invariants for the \ec es of the \sq s $\SQ^{\delta-n,5}_\lu$.

It is not difficult to see that the level curves of $I_n$ form the pencil of conics passing through four fixed points depending only on $N_5$, and the level curves of $J_n$ form the pencil of cubics passing through nine fixed points depending only on $N_5$ (one of them is on the line at infinity).  Thus generically, $\SQ^{\delta-n,5}_\lu \cong \SQ^{\delta'-n,5}_\lup$ \iff\ $(\gd)$ and $(\gdp)$ lie on the same conic and the same cubic in these pencils.  Put differently, $\EC^5_n(\gd)$ is the intersection of the conic and the cubic through $(\gd)$.  This intersection is usually six points in $(\gd)$-space, so the \ec es are usually six pairs of conjugate \sq s.

We shall describe the pencils of conics and cubics in detail in Section~\ref{Pictures}.  In particular, we shall recover the results of \cite{LO99} and \cite{Ga00} in length~5: for either $\delta$ fixed at~$0$ or the order $k = \delta-n$ fixed at~$4$, generically the only non-trivial equivalence is conjugation.

\section{Resonant results}  \label{Res Results}

In this section we answer the equivalence question in all resonant cases of length $l \le 4$, and in the self-dual resonant case of length~5.  As before, proofs are deferred to Section~\ref{Proofs}.  For the resonant modules $\D^k_\lu$, the results match those of \cite{Ga00}.

We recall relevant material from Section~\ref{Dfns}: $\SQ^{\delta-n,l}_\lu$ is resonant if its \cs\ $\{ \F_n, \F_{n+1}, \ldots, \F_{n+l-1} \}$ contains a pair of \tdm s of degrees symmetric around $\oh$.  The duality described by Lemma~\ref{trace} pairs $n$ with $2-l-n$, and hence $N_l := n + \oh l - 1$ with $-N_l$.  Thus resonance occurs for
\begin{equation*}
   n \in \bigl\{0, -\oh, -1, -{\ts\frac{3}{2}}, \ldots, 1 - \oh l \bigr\},
   \mbox{\rm\ \ \ie\ \ }
   N_l \in \bigl\{0, \pm\oh, \pm 1, \pm{\ts\frac{3}{2}}, \ldots, \pm(\oh l -1) \bigr\}.
\end{equation*}
In particular, for $l=1$ there are no resonant cases, for $l=2$ the only resonant $n$-value is~$0$, for $l=3$ the resonant $n$-values are~$-1$, $-\oh$, and~$0$, and for $l=4$ they are~$-2$, $-\frac{3}{2}$, $-1$, $-\oh$, and~$0$.

We begin by extending the definition of the \svc\ to the resonant cases and stating a proposition which resolves the majority of resonant cases of length $l \le 4$.  We will see that the self-dual cases, where $n = 1 - \oh l$ and $N_l = 0$, tend to be exceptional.

\meno {\bf Definition.}
For $n$ half-integral resonant, the simultaneous vanishing condition is the same as in the non-resonant case.  For $n$ integral resonant, two \sq s $\SQ^{\delta-n,l}_\lu$ and $\SQ^{\delta'-n,l}_\lup$ are said to satisfy the simultaneous vanishing condition if they satisfy the non-resonant \svc\ and in addition induce simultaneous vanishing of $\delta \gamma^{1/2}$.

\begin{prop} \label{resSVC}
For $l \le 3$, or $l = 4$ and $n = -\oh$ or $-\frac{3}{2}$, resonant \sq s\/ $\SQ^{\delta-n,l}_\lu$ and\/ $\SQ^{\delta'-n,l}_\lup$ are equivalent \iff\ they satisfy the \svc.
\end{prop}

Let us write the \svc\ explicitly in each of these cases.  In the self-dual case $l=2$ and $n = 0$, we have the \sq s $\SQ^{\delta,2}_\lu$ with \cs\ $\{ \F_0, \F_1 \}$.  By the proposition, two such are equivalent \iff\ they induce simultaneous vanishing of $\delta \gamma^{1/2}$.  The \ec\ where $\delta \gamma^{1/2}$ vanishes splits as $\F_0 \oplus \F_1$: the factor $\gamma^{1/2}$ reflects the splitting of self-conjugate modules into symmetric and skew-symmetric parts, and the factor $\delta$ reflects~(\ref{Dsplitting}).

For reference, we write out those functions $B_{n+i,n}(\gd)$ which arise in lengths~3 and~4, arranged in dual pairs.  Those with $i=2$ are
\begin{equation*} \ts
   B_{\frac{3}{2}, -\frac{1}{2}} = \gamma - \frac{3}{4},
\end{equation*}
\begin{equation*} \begin{array}{rclrcl} \ts
   B_{2,0} &=& \gamma - \delta - 1, \quad & \quad
   B_{1,-1} &=& \gamma + \delta - 1, \\[6pt]
   B_{\frac{5}{2}, \frac{1}{2}} &=& \gamma - 2\delta - \frac{7}{4}, \quad & \quad
   B_{\frac{1}{2}, -\frac{3}{2}} &=& \gamma + 2\delta - \frac{7}{4}, \\[6pt]
   B_{3,1} &=& \gamma - 3\delta - 3, \quad & \quad
   B_{0,-2} &=& \gamma + 3\delta - 3,
\end{array} \end{equation*}
and those with $i=3$ are
\begin{equation*} \ts
   B_{2,-1} = \gamma^{1/2} \bigl(\gamma - 3 \bigr),
\end{equation*}
\begin{equation*} \begin{array}{rclrcl} \ts
   B_{\frac{5}{2}, -\frac{1}{2}} &=& 
   \gamma^{1/2} \bigl(\gamma - \frac{3}{2} \delta - 3 \bigr), \quad & \quad
   B_{\frac{3}{2}, -\frac{3}{2}} &=& 
   \gamma^{1/2} \bigl(\gamma + \frac{3}{2} \delta - 3 \bigr), \\[6pt]
   B_{3,0} &=& \gamma^{1/2} \bigl(\gamma - 3\delta - 3 \bigr), \quad & \quad
   B_{1,-2} &=& \gamma^{1/2} \bigl(\gamma + 3\delta - 3\bigr).
\end{array} \end{equation*}

In the self-dual length~3 case $n=-\oh$, $\SQ^{\delta+\frac{1}{2},3}_\lu$ and\/ $\SQ^{\delta'+\frac{1}{2},3}_\lup$ are equivalent \iff\ they induce simultaneous vanishing of 
\begin{equation*}
   (\delta + \oh)_2\, B_{\frac{3}{2},-\frac{1}{2}}.
\end{equation*}
By \cite{Ga00, CS04}, the \ec\ where this quantity vanishes splits as $\F_{-1/2} \oplus \F_{1/2} \oplus \F_{3/2}$.  In fact, $\gamma = \frac{3}{4}$ reduces to $\lambda + \mu = \oh$ or $\frac{3}{2}$, the situation of Proposition~9.1(b) in \cite{Ga00} (see also Section~8.2 of \cite{CS04}).

In the length~3 case $n=0$, $\SQ^{\delta,3}_\lu$ and $\SQ^{\delta',3}_\lup$ are equivalent \iff\ they induce simultaneous vanishing of 
\begin{equation*} 
   \delta \gamma^{1/2}, \qquad (\delta)_2\, B_{2,0}.
\end{equation*}

In the dual case $l=3$ and $n=-1$, $\SQ^{\delta+1,3}_\lu$ and $\SQ^{\delta'+1,3}_\lup$ are equivalent \iff\ they induce simultaneous vanishing of 
\begin{equation*} 
   \delta \gamma^{1/2}, \qquad (\delta+1)_2\, B_{1,-1}.
\end{equation*}

In the length~4 resonant cases with $n = -\oh$, $\SQ^{\delta+\frac{1}{2},4}_\lu$ and\/ $\SQ^{\delta'+\frac{1}{2},4}_\lup$ are equivalent \iff\ they induce simultaneous vanishing of
\begin{equation*} \ts
   (\delta + \frac{1}{2})_2\, B_{\frac{3}{2}, -\frac{1}{2}}, \qquad
   (\delta - \frac{1}{2})_2\, B_{\frac{5}{2}, \frac{1}{2}}, \qquad
   (\delta + \frac{1}{2})_3\, B_{\frac{5}{2}, -\frac{1}{2}}.
\end{equation*}

In the dual case $l=4$ and $n=-\frac{3}{2}$, $\SQ^{\delta+\frac{3}{2},4}_\lu$ and\/ $\SQ^{\delta'+\frac{3}{2},4}_\lup$ are equivalent \iff\ they induce simultaneous vanishing of
\begin{equation*} \ts
   (\delta + \frac{1}{2})_2\, B_{\frac{3}{2}, -\frac{1}{2}}, \qquad
   (\delta + \frac{3}{2})_2\, B_{\frac{1}{2}, -\frac{3}{2}}, \qquad
   (\delta + \frac{3}{2})_3\, B_{\frac{3}{2}, -\frac{3}{2}}.
\end{equation*}

The equivalence question in length~4 at $n=0$ and $n=-2$ is resolved by the following proposition.

\begin{prop} \label{res l=4 n=0}
For $l = 4$ and $n = 0$, $\SQ^{\delta,4}_\lu$ and\/ $\SQ^{\delta',4}_\lup$ are equivalent \iff\ they induce simultaneous vanishing of
\begin{equation*} \ts
   \delta \gamma^{1/2}, \qquad
   (\delta)_2\, B_{2,0}, \qquad
   (\delta - 1)_2\, B_{3,1},
\end{equation*}
and in the case that $(\delta-1)_2\, B_{3,1} = 0$, induce in addition simultaneous vanishing of
\begin{equation*}
   (\delta)_4\, \gamma^{1/2}.
\end{equation*}

In the dual case $l=4$ and $n=-2$, $\SQ^{\delta + 2,4}_\lu$ and\/ $\SQ^{\delta' + 2,4}_\lup$ are equivalent \iff\ they induce simultaneous vanishing of
\begin{equation*} \ts
   \delta \gamma^{1/2}, \qquad
   (\delta + 1)_2\, B_{1,-1}, \qquad
   (\delta + 2)_2\, B_{0,-2},
\end{equation*}
and in the case that $(\delta+2)_2\, B_{0,-2} = 0$, induce in addition simultaneous vanishing of
\begin{equation*}
   (\delta+3)_4\, \gamma^{1/2}.
\end{equation*}
\end{prop}

Note that in the $n=0$ \dog\ case $\SQ^{3,4}_\lu = \D^3_{\lambda,\lambda+3}$, the additional condition is automatically satisfied because $\delta=3$.  This is why our results agree with those of \cite{Ga00} here.

Last we treat the self-dual cases in lengths~4 and~5.  The self-dual $l=4$ case at $n=-1$ is particularly interesting because there is a new rational invariant:
\begin{equation*}
   R(\gd)\ :=\ \frac{\gamma^{1/2}\, B_{2,-1}}{B_{2,0}\, B_{1,-1}}\ =\
   \frac{\gamma (\gamma - 3)}{(\gamma-\delta-1) (\gamma+\delta-1)}.
\end{equation*}

\begin{thm} \label{resSD4}
If\/ $\SQ^{\delta+1,4}_\lu$ and\/ $\SQ^{\delta'+1,4}_\lup$ are equivalent, then they satisfy the \svc, that is, they induce simultaneous vanishing of 
\begin{equation*}
   \delta \gamma^{1/2}, \quad
   (\delta+1)_3\, \gamma^{1/2} (\gamma-3), \quad
   (\delta)_2\, (\gamma - \delta - 1), \quad
   (\delta+1)_2\, (\gamma + \delta - 1).
\end{equation*}
If at least one of these four functions does vanish on both \sq s, then simultaneous vanishing is sufficient for equivalence.  If none of the functions vanishes, then the \sq s are equivalent \iff\ $R(\gd) = R(\gdp)$.
\end{thm}

The level curves of $R$ comprise a pencil of conics which we shall describe in Section~\ref{Pictures}.  This invariant is the reason for the exceptional behaviour of the modules $\D^3_{\lambda, \lambda+2} = \SQ^{3,4}_{\lambda,\lambda+2}$ pointed out in Theorem~3.3(4) and Section~10.1 of \cite{Ga00}: at $\delta = 2$, $R$ reduces to the monotonic function $\gamma / (\gamma+1)$, so each of these modules is equivalent only to its conjugate.

The self-dual $l=5$ case, where $n=-\frac{3}{2}$, is exceptional for the opposite reason: it is simpler than the other $l=5$ resonant cases.  We do not treat any other resonant $l \ge 5$ cases in this article, but one could do so using the simplified formulas given in Section~\ref{Proofs} for the coefficients $\o b_{m,n}$ first studied in \cite{CS04}.

\begin{prop} \label{resSD5}
The \sq s $\SQ^{\delta+\frac{3}{2},5}_\lu$ and\/ $\SQ^{\delta'+\frac{3}{2},5}_\lup$ are equivalent \iff\ they satisfy the conditions of Theorem~\ref{nonres5} with $n=-\frac{3}{2}$.
\end{prop}

\section{Lacunary \sq s}  \label{Lac Results}

In this section we answer the equivalence question in certain non-resonant lacunary cases.  Again, proofs are deferred to Section~\ref{Proofs}.  It has long been known that for $k \not = \delta$, $\Psi^k_\lu$ has a unique {\em lacunary\/} $\VR$-invariant submodule
\begin{equation*}
   \Psi^{k,\lac}_\lu, \mbox{\rm\ with \cs\ } 
   \bigl\{ \F_{\delta-k},\, \F_{\delta-k+2},\, \F_{\delta-k+3},\, \F_{\delta-k+4}, \ldots \bigr\}.
\end{equation*}
A proof may be found in \cite{Co05}, where $\Psi^{k,\lac}_\lu$ is called $\Psi^k_1(\lambda,\delta)$.  The module is given by applying the $(\lu)$-projective quantization (see Section~\ref{Proofs}) to $\F_{\delta-k} \oplus \bigoplus_{i=2}^\infty \F_{\delta-k+i}$.  The next two propositions answer the equivalence question for the simplest lacunary \sq s; they are parallel to Propositions~\ref{nonres3} and~\ref{nonres4}.

\begin{prop} \label{lac02}
The \sq s\/ $\Psi^{\delta-n,\lac}_\lu / \Psi^{\delta-n-3}_\lu$ and\/ $\Psi^{\delta'-n,\lac}_\lup / \Psi^{\delta'-n-3}_\lup$ have \cs\ $\{\F_n, \F_{n+2} \}$.  For $n\not=0$, \ie\ $N_3 \not= \oh$, they are equivalent \iff\ they induce simultaneous vanishing of~(\ref{l3poly}).

The \sq s\/ $\Psi^{\delta-n}_\lu / \Psi^{\delta-n-1,\lac}_\lu$ and\/ $\Psi^{\delta'-n}_\lup / \Psi^{\delta'-n-1,\lac}_\lup$ also have \cs\ $\{\F_n, \F_{n+2} \}$.  For $n\not= -1$, \ie\ $N_3 \not= -\oh$, they too are equivalent \iff\ they induce simultaneous vanishing of~(\ref{l3poly}).  The two cases are dual.
\end{prop}

\begin{prop} \label{lac023}
The \sq s\/ $\Psi^{\delta-n,\lac}_\lu / \Psi^{\delta-n-4}_\lu$ and\/ $\Psi^{\delta'-n,\lac}_\lup / \Psi^{\delta'-n-4}_\lup$ have \cs\ $\{\F_n, \F_{n+2}, \F_{n+3} \}$.  For $n\not\in\{-2, 0\}$, \ie\ $N_4 \not= \pm1$, they are equivalent \iff\ they induce simultaneous vanishing of
\begin{equation*}
   (\delta-N_4+1)_2\, B_{n+2,n}, \qquad
   (\delta-N_4+1)_3\, B_{n+3,n}.
\end{equation*}

The \sq s\/ $\Psi^{\delta-n}_\lu / \Psi^{\delta-n-2,\lac}_\lu$ and\/ $\Psi^{\delta'-n}_\lup / \Psi^{\delta'-n-2,\lac}_\lup$ have \cs\ $\{\F_n, \F_{n+1}, \F_{n+3} \}$, dual to the first case.  For $n\not\in\{-2, 0\}$, \ie\ $N_4 \not= \pm1$, they are equivalent \iff\ they induce simultaneous vanishing of
\begin{equation*}
   (\delta-N_4)_2\, B_{n+3,n+1}, \qquad
   (\delta-N_4+1)_3\, B_{n+3,n}.
\end{equation*}
\end{prop}

The first case in which both the numerator and the denominator of the \sq\ are lacunary has \cs\ $\{\F_n, \F_{n+2}, \F_{n+4}\}$.  Here we obtain a single rational invariant, the function $I_n$ arising in Section~\ref{l5}.

\begin{thm} \label{lac024}
The \sq s\/ $\Psi^{\delta-n,\lac}_\lu / \Psi^{\delta-n-3,\lac}_\lu$ and\/ $\Psi^{\delta'-n,\lac}_\lup / \Psi^{\delta'-n-3,\lac}_\lup$ have \cs\ $\{\F_n, \F_{n+2}, \F_{n+4}\}$.  For 
\begin{equation*} \ts
   n \not\in \bigl\{-3, -\frac{5}{2}, -\oh, 0 \bigr\}, \mbox{\rm\ \ \ie\ \ }
   N_5 \not\in \bigl\{\pm 1, \pm\frac{3}{2} \bigr\},
\end{equation*}
if they are equivalent then they induce simultaneous vanishing of
\begin{equation*} \ts
   (\delta - N_5 + \frac{3}{2})_2\, B_{n+2,n}, \quad
   (\delta - N_5 - \frac{1}{2})_2\, B_{n+4,n+2}, \quad
   (\delta - N_5 + \frac{3}{2})_4\, B_{n+4,n}.
\end{equation*}
If at least one of these three functions does vanish on both \sq s, then simultaneous vanishing is sufficient for equivalence.  If none of the functions vanishes, then the \sq s are equivalent \iff\ $I_n(\gd) = I_n(\gdp)$.
\end{thm}

Thus generically, \sq s with \cs\ $\{\F_n, \F_{n+2}, \F_{n+4}\}$ are equivalent \iff\ their parameters in $(\gd)$-space lie on the same level curve of $I_n$.  As we mentioned in Section~\ref{l5}, these level curves comprise a pencil of conics which is described in Section~\ref{Pictures}.

For completeness, we also answer the equivalence question for the two other lacunary \cs\ beginning with $\F_n$ and ending with $\F_{n+4}$:

\begin{cor} \label{lac0234}
The \sq s\/ $\Psi^{\delta-n,\lac}_\lu / \Psi^{\delta-n-5}_\lu$ and\/ $\Psi^{\delta'-n,\lac}_\lup / \Psi^{\delta'-n-5}_\lup$ have \cs\ $\{\F_n, \F_{n+2}, \F_{n+3}, \F_{n+4}\}$.  Under the same restriction on $n$ as in Theorem~\ref{lac024}, they are equivalent \iff\ they satisfy the conditions of Theorem~\ref{lac024} and in addition induce simultaneous vanishing of
$$ (\delta - N_5 + \ts\frac{3}{2})_3\, B_{n,n+3}. $$

The \sq s\/ $\Psi^{\delta-n}_\lu / \Psi^{\delta-n-3,\lac}_\lu$ and\/ $\Psi^{\delta'-n}_\lup / \Psi^{\delta'-n-5}_\lup$ have \cs\ $\{\F_n, \F_{n+1}, \F_{n+2}, \F_{n+4}\}$, dual to the last case.  Under still the same restriction on $n$, they are equivalent \iff\ they satisfy the conditions of Theorem~\ref{lac024} and in addition induce simultaneous vanishing of
$$ (\delta-N_5+\ts\frac{1}{2})_3\, B_{n+1,n+4}. $$
\end{cor}

In the case of the \cs\ $\{\F_n, \F_{n+2}, \F_{n+3}, \F_{n+5} \}$ we obtain a new rational invariant.  In order to make the symmetry of Corollary~\ref{dual} transparent, we express the relevant functions $B_{n+i,n}$ in terms of $N_6$:
\begin{equation*} \begin{array}{rcl}
   B_{n+2,n} &=& \gamma - \bigl[ (2N_6-3) \delta
   + (N_6^2 - 3N_6 + 3) \bigr], \\[6pt]
   B_{n+5,n+3} &=& \gamma - \bigl[ (2N_6+3) \delta
   + (N_6^2 + 3N_6 + 3) \bigr], \\[6pt]
   B_{n+3,n} &=& \gamma^{3/2}
   - 3\gamma^{1/2} \bigl[ (N_6 - 1) \delta + 1 \bigr], \\[6pt]
   B_{n+5,n+2} &=& \gamma^{3/2} 
   - 3\gamma^{1/2}\bigl[ (N_6 + 1) \delta + 1 \bigr].
\end{array} \end{equation*}
The invariant is
\begin{equation*}
M_n(\gd) := B_{n+5,n+2} B_{n+2,n} \big/ B_{n+5,n+3} B_{n+3,n}.
\end{equation*}

\begin{thm} \label{lac0235}
The \sq s\/ $\Psi^{\delta-n,\lac}_\lu / \Psi^{\delta-n-4,\lac}_\lu$ and\/ $\Psi^{\delta'-n,\lac}_\lup / \Psi^{\delta'-n-4,\lac}_\lup$ have \cs\ $\{\F_n, \F_{n+2}, \F_{n+3}, \F_{n+5}\}$.  For
\begin{equation*} \ts
   n \not\in \bigl\{-4, -\frac{7}{2}, -3, -2, -1, -\oh, 0 \bigr\}, \mbox{\rm\ \ \ie\ \ }
   N_6 \not\in \bigl\{0, \pm 1, \pm\frac{3}{2}, \pm 2 \bigr\},
\end{equation*}
if they are equivalent then they induce simultaneous vanishing of
\begin{eqnarray*} \ts
   (\delta - N_6 + 2)_2\, B_{n+2,n}, &\qquad&
   (\delta - N_6 - 1)_2\, B_{n+5,n+3}, \\[6pt]
   (\delta - N_6 + 2)_3\, B_{n+3,n}, &\qquad&
   (\delta - N_6)_3\, B_{n+5,n+2}.
\end{eqnarray*}
If at least one of these four functions does vanish on both \sq s, then simultaneous vanishing is sufficient for equivalence.  If none of the functions vanishes, then the \sq s are equivalent \iff\ $M_n(\gd) = M_n(\gdp)$.
\end{thm}

Note that when $\gamma\not=0$, the $\gamma^{1/2}$ factors of $M_n$ cancel, so its level curves comprise a new pencil of conics.  This pencil too will be described in Section~\ref{Pictures}.  

The analog of Corollary~\ref{lac0234} here involves both $I_n$ and $M_n$:

\begin{thm} \label{lac02345}
The \sq s\/ $\Psi^{\delta-n,\lac}_\lu / \Psi^{\delta-n-6}_\lu$ and\/ $\Psi^{\delta'-n,\lac}_\lup / \Psi^{\delta'-n-6}_\lup$ have \cs\ $\{\F_n, \F_{n+2}, \F_{n+3}, \F_{n+4}, \F_{n+5}\}$.  For
\begin{equation*} \ts
   n \not\in \bigl\{-4, -\frac{7}{2}, -3, -2, -1, -\oh, 0 \bigr\},
\end{equation*}
they are equivalent \iff\ they satisfy all of the following conditions:
\begin{enumerate}
\item[(i)]
They induce simultaneous vanishing of~(\ref{Bij}) for
\begin{equation*}
   (i,j) \in \bigl\{(2,0),\ (4,2),\ (5,3),\ (3,0),\ (5,2),\ (4,0) \bigr\}.
\end{equation*}
\smallbreak \item[(ii)]
If~(\ref{Bij}) is non-zero for $(i,j) \in \{ (2,0),\ (4,2),\ (4,0) \}$, then
\begin{equation*}
   I_n(\gd) = I_n(\gdp).
\end{equation*}
\smallbreak \item[(iii)]
If~(\ref{Bij}) is non-zero for $(i,j) \in \{ (2,0),\ (5,3),\ (3,0),\ (5,2) \}$, then
\begin{equation*}
   M_n(\gd) = M_n(\gdp).
\end{equation*}
\end{enumerate}

The \sq s\/ $\Psi^{\delta-n}_\lu / \Psi^{\delta-n-4,\lac}_\lu$ and\/ $\Psi^{\delta'-n}_\lup / \Psi^{\delta'-n-4,\lac}_\lup$ have \cs\ $\{\F_n, \F_{n+1}, \F_{n+2}, \F_{n+3}, \F_{n+5}\}$.  They are equivalent \iff\ the dual \sq s\/ $\Psi^{n+4-\delta,\lac}_\ul / \Psi^{n-2-\delta}_\ul$ and\/ $\Psi^{n+4-\delta',\lac}_\ulp / \Psi^{n-2-\delta'}_\ulp$ are equivalent.
\end{thm}

Thus generically, \sq s with \cs\ $\{\F_n, \F_{n+2}, \F_{n+3}, \F_{n+4}, \F_{n+5} \}$ are equivalent \iff\ their parameters in $(\gd)$-space lie on the same level curves of both $I_n$ and $M_n$.  Since both level curves are conics, the \ec es usually consist of four conjugate pairs of \sq s.

\section{Equivalence pencils} \label{Pictures}

In this section we examine the level curves in $(\gd)$-space of some of the rational invariants occurring in Theorems~\ref{nonres5}, \ref{resSD4}, \ref{lac024}, \ref{lac0235}, and~\ref{lac02345}.

\subsection{The invariant $R = \gamma^{1/2} B_{2,-1} / B_{2,0} B_{1,-1}$}
\begin{figure}[htb] \begin{center} \includegraphics[height=70mm]{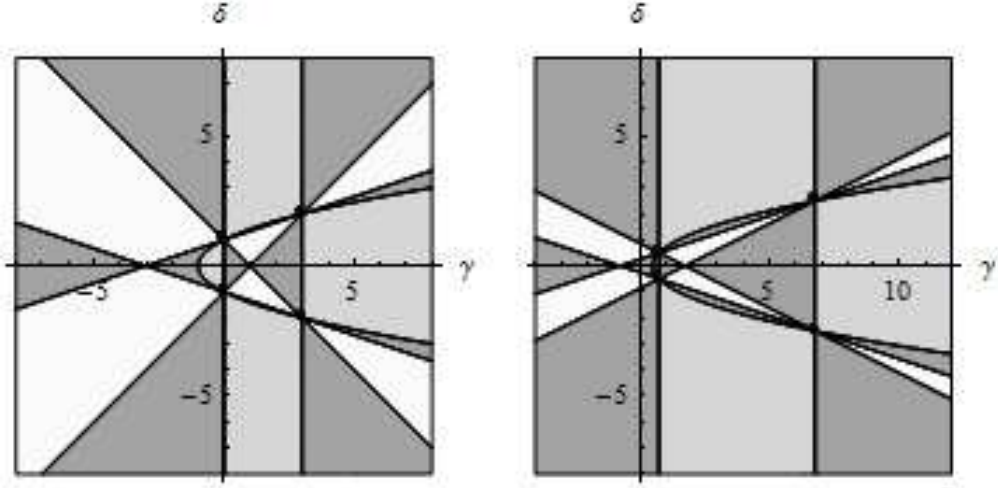} \end{center}
\caption{These are self-dual resonant invariants: at left is $R$ and at right is $I_{-3/2}$.  Depicted are the four points determining the pencil of conics and the lines and parabola bordering the different zones.  The zones of up-down and right-left hyperbolas are shaded light and dark, respectively, and the zone of ellipses is shaded grey.}  \label{R&IN0} \end{figure}

Together with the SVC, $R$ is a complete invariant for the \ec es $\EC^4_{-1}$ of length~4 \sq s with self-dual \cs\ given in Theorem~\ref{resSD4}.  It is more convenient to work with the invariant
\begin{equation*}
   \t R := \frac{R}{1 - R} = \frac{\gamma^2 - 3\gamma}{\gamma + 1 - \delta^2},
\end{equation*}
which of course is also complete.  Its level curves form the pencil of conics passing through the four points $(0, \pm 1)$ and $(3, \pm 2)$, the simultaneous zeroes of the numerator and denominator.  The conic at level $\t R$ may be written as
\begin{equation*}
   \bigl[ \gamma - \oh(\t R - 3) \bigr]^2 + \t R \delta^2 = {\ts \frac{1}{4}} (\t R + 1)(\t R + 9).
\end{equation*}
For $\t R$ real, the real points form ellipses for $\t R > 0$, hyperbolas for $\t R < 0$, a parabola for $\t R = \infty$, two parallel lines for $\t R = 0$, and two crossing lines for $\t R = -1$ or~$-9$.  These different zones are delineated by shadings in Figure~\ref{R&IN0}.

\subsection{The invariant $I_n = B_{n+4,n} / B_{n+4,n+2} B_{n+2,n}$}
This is the simpler of the two continuous invariants for the \ec es $\EC^5_n$ of length~5 \sq s treated in Theorem~\ref{nonres5} and Proposition~\ref{resSD5}, and together with the SVC it is a complete invariant for the \ec es of lacunary modules treated in Theorem~\ref{lac024}.  Like $R$, its level curves form a pencil of conics.  To describe them it will be convenient to use the coordinates $(\tgfd)$, where $\t\gamma_5 := \gamma - 2 N_5 \delta$, as the conics are all in standard orientation in these coordinates.  Some computation gives
\begin{align*}
   &B_{n+3 \pm 1, n+1 \pm 1}\ =\ 
   \bigl[ \t\gamma_5 - N_5^2 - {\ts\frac{7}{4}} \bigr]
   \mp 2 \bigl[ \delta + N_5 \bigr], \\[6pt]
   &B_{n+4,n}\ =\ 
   \bigl[ \t\gamma_5 + N_5^2 - {\ts\frac{15}{4}} \bigr]^2
   -4 \bigl[ N_5 \delta - {\ts\frac{1}{5}} N_5^2 + {\ts\frac{6}{5}} \bigr]^2
   - {\ts\frac{9}{25}} (N_5^2 - 1) (4N_5^2 - 9).
\end{align*}

The pencil of conics is determined by the four simultaneous zeroes of $B_{n+4,n}$ and $B_{n+4,n+2} B_{n+2,n}$.  In $(\tgfd)$-coordinates, these are
\begin{equation*}
   \bigl( N_5^2 \pm 4N_5 + {\ts\frac{27}{4}},\ N_5 \pm {\ts\frac{5}{2}} \bigr), \qquad
   \bigl( N_5^2 \pm {\ts\frac{4}{5}} N_5 + {\ts\frac{3}{4}},\
   -{\ts\frac{3}{5}} N_5 \mp \oh \bigr),
\end{equation*}
a quadrilateral with center $\bigl( N_5^2 + {\ts\frac{15}{4}},\ {\ts\frac{1}{5}} N_5 \bigr)$.  The slopes of the six lines connecting the vertices are $\pm 2$, $\pm 3$, and $\pm \frac{8}{5} N_5$, from which it follows that the quadrilateral is cyclic (inscribed).  Note that at $N_5 = 0$ and as $N_5 \to \infty$ it becomes a trapezoid.  When $N_5$ is $\pm 5/4$ or $\pm 15/8$ two of the four points coincide; these are the only values of $N_5$ for which this occurs.

We remark that since $B_{n+4,n}$ is a conic, it is at first surprising that the four vertices are all rational in $N_5$.  However, this is explained by the final paragraphs of \cite{Co05}: one of the two simultaneous zeroes of $B_{n+2,n}$ and $B_{n+4,n}$ is also a zero of $B_{n+3,n} / \gamma^{1/2}$, a linear function, so it is rational, forcing the other to be rational.  The rationality of the simultaneous zeroes of $B_{n+4,n+2}$ and $B_{n+4,n}$ follows by duality.

Completing the square in $\t\gamma_5$ and $\delta$, a long calculation gives the following form of the conic at level $I_n$:
\begin{align}
   & (I_n - N_5^2) \bigl[ (I_n-1) (\t\gamma_5 - N_5^2 - {\ts\frac{15}{4}})
   + 2(I_n - N_5^2) \bigr]^2 \nonumber \\[6pt]
   & - 4(I_n - 1) \bigl[ (I_n - N_5^2) (\delta - {\ts\frac{1}{5}} N_5)
   + {\ts\frac{6}{5}} (I_n - 1) N_5 \bigr]^2 \label{Ipencil} \\[6pt]
   =\ & 
   -{\ts\frac{1}{25}} (N_5^2 - 1)
   \bigl[ 9(I_n - 1) - 4(I_n - N_5^2) \bigr]
   \bigl[ 16 N_5^2 (I_n - 1) - 25 (I_n - N_5^2) \bigr]. \nonumber
\end{align}
When $I_n$ is~$1$ or~$N_5^2$, taking limits in the obvious way gives a parabola.  We obtain ellipses when $(I_n - 1) (I_n - N_5^2)$ is negative and hyperbolas when it is positive.  When $I_n$ is~$\infty$, $-\frac{4}{5}(N_5^2 - \frac{9}{4})$, or $-9N_5^2 / (16N_5^2 - 25)$, the right side is zero and we obtain a degenerate hyperbola: either a pair of opposite sides of the quadrilateral or its diagonal.

\begin{figure}[htb] \begin{center} \includegraphics[height=70mm]{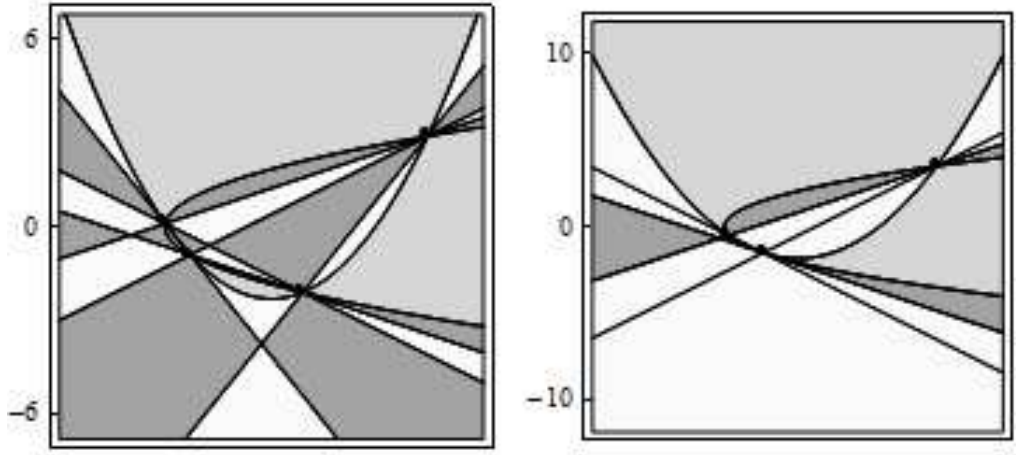} \end{center}
\caption{Zones of ellipses and hyperbolas for $I_{-1}$ (left) and $I_{-1/4}$ (right) in $(\t\gamma_5, \delta)$-coordinates.  The zones are shaded as in Figure~\ref{R&IN0}.}  \label{IN.5&IN1.25} \end{figure}

These zones may be seen for some values of $n$ in Figures~\ref{R&IN0}, \ref{IN.5&IN1.25}, and~\ref{IN12&MN6.35}.  The zone of vertical hyperbolas is light, the zone of horizontal hyperbolas is dark, and the zone of ellipses is grey.  Some simple rules describe the zones: crossing either of the two parabolas toggles between the elliptical zone and the hyperbolic zones, and crossing a line toggles between the two hyperbolic zones.  Crossing over a parabola out of a hyperbolic zone and then back over the same parabola leads to the same hyperbolic zone, whereas crossing back over the other parabola leads to the other hyperbolic zone.

On the right in Figure~\ref{R&IN0} is the self-dual case $I_{-3/2}$ relevant to Proposition~\ref{resSD5}, where $N_5 = 0$.  In Figure~\ref{IN.5&IN1.25} we show the resonant case $I_{-1}$ occurring in Theorem~\ref{lac024}, where $N_5 = 1/2$, and the case $I_{-1/4}$, where $N_5 = 5/4$ and the quadrilateral has a double vertex, causing two of the degenerate hyperbolas to coincide.  On the left in Figure~\ref{IN12&MN6.35} is $I_{10.5}$, where one can see the horizontal parabola approaching a pair of lines.  In all cases the magnification is the same on both axes, but in the non-self-dual cases we have not numbered the $\t\gamma_5$ axis because it is shifted in order to center on the quadrilateral.  Observe that the rules describing the zones sometimes appear to be violated, because some of the zones are too thin to be seen in the figures.  For example, in the picture of $I_{-1}$ the line and the vertical parabola connecting the two left points cannot be distinguished, so the light zone between them is invisible.

Recall from Theorem~\ref{lac024} that $I_n$ is meaningless when $N_5$ is one of the resonant values $\pm 1$ or $\pm \frac{3}{2}$.  The values $N_5 = \pm 1$ are particularly special: there $B_{n+4,n}$ and $B_{n+4,n+2} B_{n+2,n}$ are equal by Corollary~\ref{B factors}, so $I_n$ is identically~$1$.  It follows that $B_{n+4,n} - B_{n+4,n+2} B_{n+2,n}$ is divisible by $N_5^2 - 1$, so we define
\begin{equation*}
   B_{420}\ :=\ \frac{B_{n+4,n} - B_{n+4,n+2} B_{n+2,n}}{4 (N_5^2 - 1)}
   \ =\ \gamma - \bigl[ \delta^2 + {\ts\frac{8}{5}} N_5 \delta
   + {\ts\frac{2}{5}} N_5^2 + \oh \bigr].
\end{equation*}
Just as we replaced the invariant $R$ by $\t R$, we can replace $I_n$ by
\begin{equation*}
   \t I_n := \frac{4 (N_5^2 - 1)}{I_n - 1}
   = \frac{B_{n+4,n+2} B_{n+2,n}}{B_{420}}.
\end{equation*}
Where $I_n$ alone is concerned, this simplifies some calculations, but not the final result~(\ref{Ipencil}).  However, it is a useful idea when considering $I_n$ together with $J_n$, which we now do.

\subsection{The invariants $I_n$ and $J_n = B_{n+4,n} B_{n+3,n+1} / B_{n+4,n+1} B_{n+3,n}$}
By Theorem~\ref{nonres5} and Proposition~\ref{resSD5}, these two invariants together with the SVC completely classify the \ec es $\EC^5_n$ unless $N_5$ takes on one of the resonant values $\pm\oh$, $\pm 1$, and $\pm\frac{3}{2}$.  Applying Corollary~\ref{B factors}, we find that $J_n$ is identically~$1$ for $N_5 = \pm\frac{3}{2}$; for example, at $N_5 = \frac{3}{2}$ we have $n=0$, $B_{4,0} = \gamma^{1/2} B_{4,1}$, and $B_{3,0} = \gamma^{1/2} B_{3,1}$.  Hence the difference $B_{n+4,n} B_{n+3,n+1} - B_{n+4,n+1} B_{n+3,n}$ is divisible by $N_5^2 - \frac{9}{4}$.  We define
\begin{align*}
   & B_{4310}\ :=\ \frac{B_{n+4,n} B_{n+3,n+1}
   - B_{n+4,n+1} B_{n+3,n}}{N_5^2 - \frac{9}{4}} \\[6pt]
   &\ =\ \gamma^2 - \gamma \bigl[ \delta^2 + {\ts\frac{12}{5}} N_5 \delta
   + {\ts\frac{13}{5}} N_5^2 + {\ts\frac{3}{4}} \bigr]
   + {\ts\frac{3}{5}} \bigl[ 2 N_5 \delta + N_5^2  + {\ts\frac{3}{4}} \bigr]
   \bigl[ 4 N_5 \delta + N_5^2 + {\ts\frac{15}{4}}\bigr].
\end{align*}

Now combine $B_{420}$ and $B_{4310}$ in the following two ways: set
\begin{align*}
   & B^-_{43210}\ :=\ 5 (B_{4310} - B_{n+3,n+1} B_{420}) \\[6pt]
   &\ =\  2 \gamma \bigl[ 3 N_5 \delta - 3 N_5^2 + {\ts\frac{5}{4}} \bigr]
   - \bigl[ 2 N_5 \delta + N_5^2  + {\ts\frac{3}{4}} \bigr]
   \bigl[ 5 \delta^2 - 4 N_5 \delta - N_5^2 - {\ts\frac{35}{4}}\bigr], \\[8pt]
   & B^+_{43210}\ := \oh (B_{4310} + B_{n+3,n+1} B_{420}) \\[6pt]
   &\ =\ \gamma^2 - \gamma \bigl[ \delta^2 + 3 N_5 \delta + 2 N_5^2 + 1 \bigr]
   + \oh \bigl[ 2 N_5 \delta + N_5^2 + {\ts\frac{3}{4}} \bigr]
   \bigl[ \delta^2 + 4 N_5 \delta + N_5^2 + {\ts\frac{11}{4}} \bigr].
\end{align*}
We will prove the following proposition at the end of Section~\ref{Res Proofs}.  Define
\begin{equation*}
   \t J_n\ :=\ B^+_{43210} / B^-_{43210}.
\end{equation*}

\begin{prop} \label{general invnt}
For any linearly independent elements $x$ and $y$ of\/ $\bC^3$,
\begin{equation*}
   \frac{x_1 B_{n+4,n+1} B_{n+3,n} + x_2 B_{n+4,n} B_{n+3,n+1}
   + x_3 B_{n+4,n+2} B_{n+3,n+1} B_{n+2,n}}
   {y_1 B_{n+4,n+1} B_{n+3,n} + y_2 B_{n+4,n} B_{n+3,n+1}
   + y_3 B_{n+4,n+2} B_{n+3,n+1} B_{n+2,n}}
\end{equation*}
is invariant under equivalence.  In particular, $\t J_n$ is an invariant.  In the situation of Part~(v) of Theorem~\ref{nonres5}, $\t I_n$ and $\t J_n$ form a complete set of invariants.
\end{prop}

We write $\t I_n$ and $\t J_n$ together for reference:
\begin{align*}
   \t I_n &= \frac{\gamma^2
   - 2\gamma \bigl[ 2 N_5 \delta + N_5^2 + {\ts\frac{7}{4}} \bigr]
   + \bigl[ (2 N_5 \delta + N_5^2 + {\ts\frac{7}{4}})^2
   - 4(\delta + N_5)^2 \bigr]}
   {\gamma - \bigl[ \delta^2 + {\ts\frac{8}{5}} N_5 \delta
   + {\ts\frac{2}{5}} N_5^2 + \oh \bigr]}, \\[6pt]
   \t J_n &= \frac{\gamma^2
   - \gamma \bigl[ \delta^2 + 3 N_5 \delta + 2 N_5^2 + 1 \bigr]
   + \oh \bigl[ 2 N_5 \delta + N_5^2 + {\ts\frac{3}{4}} \bigr]
   \bigl[ \delta^2 + 4 N_5 \delta + N_5^2 + {\ts\frac{11}{4}} \bigr]}
   {2 \gamma \bigl[ 3 N_5 \delta - 3 N_5^2 + {\ts\frac{5}{4}} \bigr]
   - \bigl[ 2 N_5 \delta + N_5^2  + {\ts\frac{3}{4}} \bigr]
   \bigl[ 5 \delta^2 - 4 N_5 \delta - N_5^2 - {\ts\frac{35}{4}}\bigr]}.
\end{align*}

Let us remark that eliminating $\gamma$ from these equations generically yields a sextic in $\delta$ with coefficients depending on $N_5$, $\t I_n$, and $\t J_n$: simply clear denominators in both equations, take the difference to eliminate $\gamma^2$, solve for $\gamma$ in terms of $\delta$, plug the result into the formula for $\t I_n$, and clear denominators again.  As we observed at the end of Section~\ref{Nonres Results}, this is expected: each equivalence class $\EC^5_n(\gd)$ is the intersection of a conic and a cubic, so generic classes contain six points.

We have not depicted the family of cubics involved because it is not unique.  The family of conics is necessarily the set of level curves of $I_n$, which is the same as the set of level curves of $\t I_n$, but the family of cubics can be altered by adding multiples of $I_n$ to the cubic invariant.  It is always a pencil determined by nine points, but there is a 1-parameter family of choices for this set of points and we did not find any ``best'' choice.  For example, $J_n$, $K_n$, and $\t J_n$ all yield different choices.

We now compare our results in length~5 with those of \cite{LO99} and \cite{Ga00}.  Recall that \cite{LO99} treats only those \sq s of \psdog s with $\lambda = \mu$, \ie\ $\delta=0$.  In fact a further restriction is imposed: only \sq s of real order $k$ are admitted, so $n$ and hence $N_5$ must be real.  The following result is stated: for $k$ real and not non-negative half-integral, the \sq s $\SQ^{k,5}_\ll$ and $\SQ^{k,5}_{\lambda',\lambda'}$ are equivalent \iff\ they are either equal or conjugate, \ie\ $\lambda'$ is either $\lambda$ or $1 - \lambda$.

On the other hand, as we have seen, \cite{Ga00} allows $\lambda$ and $\mu$ to vary independently but treats only \dog\ modules.  In fact, both $\lambda$ and $\mu$ are required to be real, so again, $n$ and $N_5$ must be real.  The following result is stated: for $\lambda$ and $\mu$ real, the \dog\ modules $\D^4_\lu$ and $\D^4_{\lup}$ are equivalent \iff\ they are either equal or conjugate.

Let us use Theorem~\ref{nonres5} to generalize these results in the nonresonant cases.  We begin with the \cite{Ga00} result, where the analysis is simpler.  The statement is true for all complex values of $\lambda$ and $\mu$.

\begin{prop} \label{check Ga00}
{\em (\cite{Ga00})}  Two non-resonant modules $\D^4_\lu$ and $\D^4_{\lup}$ are equivalent \iff\ they are either equal or conjugate, \ie\ $(\lup)$ is either $(\lu)$ or $(1-\mu, 1-\lambda)$.
\end{prop}

\meno {\it Proof.\/}
The order $k = \delta - n$ is~$4$ at $\delta = N_5 + \frac{5}{2}$.  Evaluating $\t J_n$ here (we used a software package) gives
\begin{equation*}
   \frac{ 2\gamma - (N_5 + 1)(6N_5 + 1) }{ 5(6N_5 + 1) }.
\end{equation*}
Thus at this value of $\delta$, the $\gamma$-linear denominator of the original expression for $\t J_n$ divides its $\gamma$-quadratic numerator for all values of $N_5$.  (Probably there is a conceptual explanation of this.)  Thus for $N_5 \not= -\frac{1}{6}$, the invariant $\t J_n$ determines $\gamma$, proving the result.  At $N_5 = -\frac{1}{6}$, the alternate invariant $B_{4310} / B_{n+4,n+1} B_{n+3,n}$ reduces to $\gamma + \frac{5}{3}$, yielding a similar proof.  $\Box$

\medbreak
In the case of the \cite{LO99} result, our result is slightly different.  We find that for most values of the \psdog\ order $k$, exactly one equivalence class $\EC^5_n$ contains two points of the form $(\gamma,0)$, and all the others contain only one.  Thus if we fix $\delta$ at zero, then for most choices of the \cs\ all but one of the \sq s is equivalent only to its conjugate, but one equivalence class consists of two conjugate pairs of \sq s.

To state the result concisely we make some preliminary definitions.  Set
\begin{align*}
   a_I &:= -2(N_5^2 + \ts\frac{7}{4}),  \quad&\quad
   b_I &:= (N_5^2 + 2N_5 + \ts\frac{7}{4}) (N_5^2 - 2N_5 + \ts\frac{7}{4}), \\[6pt]
   c_I &:= 1, \quad&\quad
   d_I &:= -\ts\frac{2}{5} (N_5^2 + \ts\frac{5}{4}) \\[6pt]
   a_J &:= -2(N_5^2 + \oh), \quad&\quad
   b_J &:= \oh (N_5^2 + \ts\frac{3}{4}) (N_5^2 + \ts\frac{11}{4}), \\[6pt]
   c_J &:= -6 (N_5^2 - \ts\frac{5}{12}), \quad&\quad
   d_J &:= (N_5^2 + \ts\frac{3}{4}) (N_5^2 + \ts\frac{35}{4}),
\end{align*}
so that at $\delta = 0$ we have
\begin{equation*}
   \t I_n = \frac{\gamma^2 + a_I \gamma + b_I}{c_I \gamma + d_I}, \qquad
   \t J_n = \frac{\gamma^2 + a_J \gamma + b_J}{c_J \gamma + d_J}.
\end{equation*}
Define polynomials $f_+(N_5)$, $f_-(N_5)$, and $E(N_5)$ by
\begin{align*}
   f_+ &:= \ts -\frac{14}{5} N_5^6 + \frac{38}{5} N_5^4 +
   \frac{145}{8} N_5^2 + \frac{175}{16}, \\[6pt]
   f_- &:= \ts\frac{994}{25} N_5^{10} - \frac{5549}{20} N_5^8 + 322\, N_5^6
   + \frac{29025}{32} N_5^4 - \frac{12125}{128} N_5^2 - \frac{380625}{1024}, \\[6pt]
   E &:= c_I d_J - c_J d_I\, =\, -\ts\frac{7}{5} (N_5^2 - \frac{25}{4}) (N_5^2 + \frac{25}{28}).
\end{align*}

\begin{prop}
We assume that $N_5$ is not one of the resonant values $\pm\oh$, $\pm 1$, $\pm\frac{3}{2}$, and also that we are in the situation of Part~(v) of Theorem~\ref{nonres5}.  If $N_5$ is either a root of $f_-$ or any of the values
\begin{equation*}
   \pm\frac{5}{2}, \quad \pm\frac{5}{6}, \quad \pm\frac{i \sqrt{3}}{2}, \quad
   \pm\frac{5 i}{2\sqrt{7}},
\end{equation*}
then each \ec\/ $\EC^5_n$ contains a unique element $(\gamma,0)$ with $\delta$-coordinate zero.

At all other values of $N_5$, each \ec\/ $\EC^5_n$ contains a unique element $(\gamma,0)$ with $\delta$-coordinate zero with one exception: $(\gamma, 0)$ and $(\gamma', 0)$ are in the same \ec\ for the two values $\gamma$ and $\gamma'$ determined by the equations
\begin{equation*}
   \gamma + \gamma' = f_+/E, \qquad (\gamma - \gamma')^2 = f_-/E^2.
\end{equation*}
\end{prop}

\meno {\it Proof.\/}
Use software to check that $\t I_n$ reduces to a linear polynomial in $\gamma$ \iff\ $N_5^2$ is $\frac{25}{4}$ or $\frac{25}{36}$, and the same occurs for $\t J_n$ \iff\ $N_5^2$ is $\frac{9}{4}$, $\frac{25}{4}$, $-\frac{3}{4}$, or $-\frac{25}{28}$.  At these values the proof goes as for Proposition~\ref{check Ga00}.  

Now suppose that $N_5$ does not take on any of these values, and note that this implies $E \not= 0$.  Under our assumptions, two values $\gamma$ and $\gamma'$ give equivalent $\delta=0$ \sq s \iff\ they give equal values of $\t I_n$ and $\t J_n$.  Write this condition in terms of $\gamma_\pm := \gamma \pm \gamma'$, clear denominators, assume that $\gamma_- \not=0$, and divide by it.  One obtains the equations in the proposition.  The reader may check that the results are the same even at $N_5^2 = \frac{5}{12}$, where $c_J = 0$ so $\t J_n$ has $\gamma$-constant denominator.  $\Box$

\medbreak
We remark that the condition that the situation of Part~(v) of Theorem~\ref{nonres5} obtains could be sharpened to the condition that $N_5$ not be a root of any of a collection of non-trivial polynomials.  We did not investigate the existence of double \ec es when $N_5$ is a root of one of these polynomials.  It is easy to use our methods to analyze the situation for any given value of $N_5$, but analyzing all cases would be arduous and as far as we can see not interesting.

\subsection{The invariant $M_n = B_{n+5,n+2} B_{n+2,n} / B_{n+5,n+3} B_{n+3,n}$}
Together with the SVC, this is a complete invariant for the lacunary equivalence classes treated in Theorem~\ref{lac0235}.  Its level curves form a pencil of conics which we will only discuss briefly.  Corollary~\ref{B factors} implies that at $N_6 = 0$, $B_{3,0} B_{0,-2} = B_{3,1} B_{1,-2}$, so we define
\begin{equation*}
   B_{5320}\ :=\ {\ts\frac{1}{6}} (B_{n+5,n+2} B_{n+2,n} - B_{n+5,n+3} B_{n+3,n}) / N_6
   \ =\ \gamma^{3/2} - \gamma^{1/2} \bigl[ \delta^2 + 2 N_6 \delta + 3 \bigr].
\end{equation*}

Proceeding along the same lines as before, we replace $M_n$ by the invariant
\begin{equation*}
   \t M_n\ :=\ \frac{6 N_6 M_n}{M_n - 1}\ =\ \frac{B_{n+5,n+2} B_{n+2,n}}{B_{5320}}.
\end{equation*}
Setting $\t\gamma_6 := \gamma - \frac{5}{2} N_6 \delta$, we find that the conics are all in standard orientation in the coordinates $(\t\gamma_6, \delta)$, as
\begin{equation*}
   \t M_n\ :=\ \frac{\bigl[ \t\gamma_6 - \bigl( \oh N_6 + 3 \bigr) \delta - 3 \bigr]
   \bigl[ \t\gamma_6 + \bigl( \oh N_6 + 3 \bigr) \delta - (N_6^2 - 3N_6 + 3) \bigr]}
   {\t\gamma_6 - \delta^2 + \oh N_6 \delta - 3}.
\end{equation*}

\begin{figure}[hbt] \begin{center} \includegraphics[height=70mm]{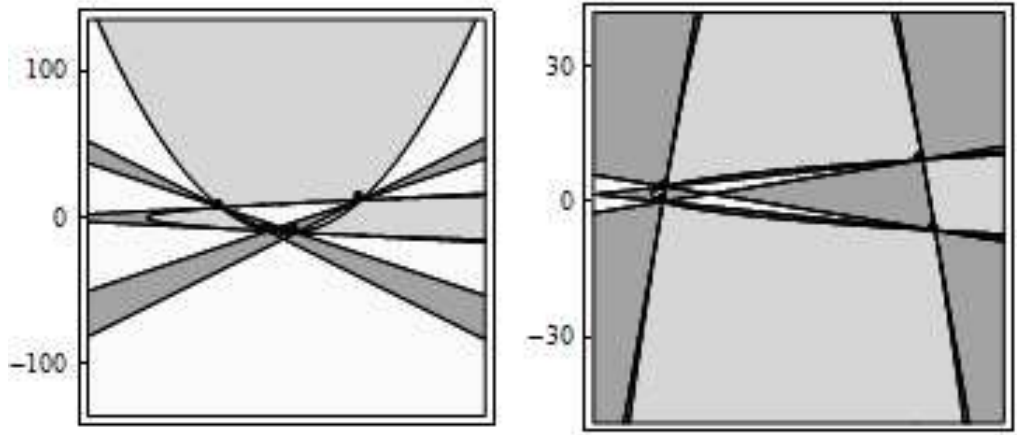} \end{center}
\caption{Zones of ellipses and hyperbolas for $I_{10.5}$ in $(\t\gamma_5, \delta)$-coordinates (left) and $M_{4.35}$ in $(\t\gamma_6, \delta)$-coordinates (right).  The zones are shaded as in Figure~\ref{R&IN0}.}  \label{IN12&MN6.35} \end{figure}

Just as for $I_n$, the four simultaneous zeroes of the numerator and denominator determine the pencil of conics.  In $(\t\gamma_6, \delta)$-coordinates, these zeroes are
\begin{equation*}
   \bigl( 3, 0 \bigr), \qquad
   \bigl( -N_6^2 + 3, -N_6 \bigr), \qquad
   \bigl( \oh N_6^2 \pm {\ts\frac{9}{2}} N_6 + 12, N_6 \pm 3 \bigl).
\end{equation*}
The slopes of the six lines determined by this quadrilateral are $\pm 2/3$ and $\pm 2 / (N_6 \pm + 6)$.  Thus the quadrilateral is cyclic and becomes a trapezoid at $N_6 = \pm 6$ and as $N_6 \to \infty$.  Two of its vertices coincide when $N$ is $0$, $\pm \frac{3}{2}$, or $\pm 3$; the coincident vertices are at either $(3,0)$ or $(\frac{51}{8}, \pm\frac{3}{2})$ in all cases.

The level curve of $\t M_6$ is a horizontal parabola at $\t M_6 = \infty$ and a vertical parabola at $\t M_6 = (\oh N_6 + 3)^2$.  At $N_6 = \pm 6$, where $n = 4$ and $-8$, the vertical parabola degenerates into two parallel lines because the quadrilateral is a trapezoid.  At these values of $N_6$ the direction of the vertical parabola reverses, which is why it is upside down on the right side of Figure~\ref{IN12&MN6.35}.

\section{Proofs} \label{Proofs}

\subsection{The non-resonant case} \label{Res Proofs}
Recall from Section~\ref{Dfns} that for $n$ non-resonant, the eigenvalues of the Casimir operator $Q$ of the projective subalgebra $\da$ on the \cs\ $\{ \F_n, \ldots, \F_{n+l-1} \}$ of $\SQ^{\delta-n,l}_\lu$ are distinct.  In this case one has the {\em projective quantization,\/} the unique $\da$-equivalence
\begin{equation} \label{PQ}
   \PQ_\lu: \bigoplus_{i=0}^{l-1} \F_{n+i} \to \SQ^{\delta-n,l}_\lu
\end{equation}
which {\em preserves symbols,\/} in the sense that it maps $\F_{n+i}$ to operators of order $\delta-n-i$ and composing it with the natural symbol map $dx^\delta f \px^{\delta-n-i} \mapsto dx^{n+i} f$ gives the identity on $\F_{n+i}$.

This map has been the subject of many articles, notably \cite{CMZ97} and \cite{LO99}.  It can be used to write the action $L_\lu$ of $\VR$ on $\SQ^{\delta-n,l}_\lu$ in an explicitly $\da$-diagonal manner:

\meno {\bf Definition.}
Let $\pi^\lu$ be the \r\ of $\VR$ on $\bigoplus_{i=0}^{l-1} \F_{n+i}$ given by
\begin{equation} \label{pilu}
   \pi^\lu(X) := \PQ_\lu^{-1} \circ \big(L_\lu(X)\big|_{\SQ^{\delta-n,l}_\lu}\big) \circ \PQ_\lu.
\end{equation}
Regard $\pi^\lu$ as an $l \times l$ matrix with entries
\begin{equation*}
   \pi^\lu_{n+i,n+j}: \VR \to \Hom(\F_{n+j}, \F_{n+i}).
\end{equation*}

\medbreak
Recall that if $W$ is any $\VR$-module, a map $\beta:\VR \to W$ is an {\em $\da$-relative $W$-valued 1-cochain\/} if it is $\da$-covariant and zero on $\da$.  The space of all such maps is denoted by $C^1\bigl(\VR, \da, W \bigr)$.

As noted in \cite{LO99}, the invariance of the order filtration on $\Psi^{\delta-n}_\lu$ implies that $\pi^\lu$ is lower triangular with tensor density actions on the diagonal, and $\PQ$'s $\da$-covariance forces the subdiagonal entries of $\pi^\lu$ to be $\da$-relative $\Hom(\F_{n+j}, \F_{n+i})$-valued 1-cochains:

\begin{lemma} \label{lower tri}
\begin{enumerate}
\item[(i)]
For $i<j$, $\pi^\lu_{n+i,n+j} = 0$.
\smallbreak \item[(ii)]
At $j=i$, $\pi^\lu_{n+i,n+i} = L_{n+i}$.
\smallbreak \item[(iii)]
For $i>j$, $\pi^\lu_{n+i,n+j} \in C^1\bigl(\VR, \da, \Hom(\F_{n+j}, \F_{n+i}) \bigr)$.
\end{enumerate}
\end{lemma}

Our proofs hinge on explicit formulas for the subdiagonal entries $\pi^\lu_{n+i,n+j}$ for $i-j \le 4$.  These formulas are given in Section~7.8 of \cite{LO99}, and they can be deduced for all $i>j$ from the results of \cite{CMZ97}.  We now state them in a form similar to but more useful than the form in which they are stated in \cite{Co05}.  The next definition, lemma, and corollary go back essentially to \cite{LO99}; they may be found in roughly the current notation in \cite{CS04} and \cite{Co05}.  

\meno {\bf Definition.}
For $\delta \in 2+\bN$, let $\beta_\lu: \VR \to \D_\lu^{\delta-2}$ be
\begin{equation*}
   \beta_\lu(g\px) := \PQ_\lu(dx^2 g''').
\end{equation*}

\begin{lemma} \label{beta props}
\begin{enumerate}
\item[(i)]
$\beta_\lu$ is $\da$-relative, and for $r\ge 2$,
\begin{equation*}
   \beta_\lu(x^{r+1}\px) := 6(r-2)!^{-1} L_\lu (x^2 \px)^{r-2} (dx^\delta \px^{\delta-2}).
\end{equation*}
\smallbreak \item[(ii)]
For $\delta \not\in 2+\bN$, $C^1\bigl(\VR, \da, \Hom(\F_\lambda, \F_\mu) \bigr) = 0$.
\smallbreak \item[(iii)]
For $\delta \in 2+\bN$, $C^1\bigl(\VR, \da, \Hom(\F_\lambda, \F_\mu) \bigr) = \bC \beta_\lu$.
\smallbreak \item[(iv)]
Under conjugation, $\beta_{\lu}^* = (-1)^\delta \beta_{1-\mu,1-\lambda}$.
\end{enumerate}
\end{lemma}

\begin{cor} \label{entries}
\begin{enumerate}
\item[(i)]
$\pi^\lu_{n+j+1,n+j} = 0$.
\smallbreak \item[(ii)]
For $i \ge j+2$, there are scalars $b_{n+i,n+j}(\lu)$ such that
\begin{equation} \label{pi b beta}
   \pi^\lu_{n+i,n+j} = b_{n+i,n+j}(\lu)\, \beta_{n+j,n+i}.
\end{equation}
\end{enumerate}
\end{cor}

In order to give these scalars in the most useful form, we define intermediate scalars $B_{m+r,m}$ for arbitrary $m \in \bC$ and $r \in 2+\bN$:
\begin{eqnarray*}
   B_{m+r, m} &:=& 3^{r/2}\, (\lambda+\mu-1-m)_r \\[6pt]
   && +\ 3^{r/2}\, \frac{2m-1+r}{r^2-1}\, \sum_{s=0}^{r-2}
   {r+1 \choose s}\, (2m-3+r)_{r-s-2}\, (\lambda + \mu-1-m)_s \\[6pt]
   && \times\ \bigl[ (r-s+1)(\lambda+\mu-1) - (r-s-1)\delta - (2m+r+s-1) \bigr].
\end{eqnarray*}
The reader may check that this formula generalizes~(\ref{3Bs}).  The point of the common factor $3^{r/2}$ is to make the expression monic when expressed in terms of $\gamma$.

\begin{thm} \label{CMZ} \cite{CMZ97}
The scalars $b_{m+r,m}(\lu)$ are given by
\begin{equation*}
   b_{m+r,m}\ =\ \frac{3^{-r/2}\, (-1)^{r-1}\, (r^2-1)\, (\delta-m)_r\, B_{m+r,m}}
   {12\, (2m-2+2r)_{r-2}\, (2m-1+r)\, (2m-3+r)_{r-2}}\,.
\end{equation*}
\end{thm}

The next proposition gives the parity of the $b_{m+r,m}$ in $\gamma^{1/2}$, and in simultaneously reversing the sign of $\delta$ and reflecting across the {\em antidiagonal\/} $m+r=1-m$.  It is a manifestation of Lemmas~\ref{conj} and~\ref{trace} and was observed in \cite{CMZ97}, except for the last sentence, which is obvious by continuity.  We remark that these parities are not easy to prove directly.

\begin{prop} \label{parity}
\begin{enumerate}
\smallbreak \item[(i)]
$\pi_{m+r,m}(1-\mu, 1-\lambda) = (-1)^r \pi_{m+r,m}(\lu)$.
\smallbreak \item[(ii)]
Under conjugation, $\pi_{1-m,1-m-r}(\ul) = -\pi_{m+r,m}(\lu)^*$.
\smallbreak \item[(iii)]
As functions of $(\gamma^{1/2}, \delta)$, both $b_{m+r,m}$ and $B_{m+r,m}$ are of parity~$r$:
\begin{eqnarray*}
   b_{m+r,m}(-\gamma^{1/2}, \delta) &=&
   (-1)^r b_{m+r,m}(\gamma^{1/2}, \delta), \\[6pt]
   B_{m+r,m}(-\gamma^{1/2}, \delta) &=&
   (-1)^r B_{m+r,m}(\gamma^{1/2}, \delta).
\end{eqnarray*}
\smallbreak \item[(iv)]
Under simultaneously reversing the sign of $\delta$ and reflecting across the antidiagonal, $b$ is of parity $r-1$ and $B$ is even:
\begin{eqnarray*}
   b_{1-m,1-m-r}(\gamma^{1/2}, -\delta) &=&
   (-1)^{r-1} b_{m+r,m}(\gamma^{1/2}, \delta), \\[6pt]
   B_{1-m,1-m-r}(\gamma^{1/2}, -\delta) &=&
   B_{m+r,m}(\gamma^{1/2}, \delta).
\end{eqnarray*}
\end{enumerate}
These parities hold for $B_{\bullet,\bullet}$ even when $b_{\bullet,\bullet}$ is undefined.
\end{prop}

We shall usually write $B_{m+r,m}$ as a function of $(\gd)$, it being understood that in fact it is a function of $(\gamma^{1/2},\delta)$.  We now use Corollary~\ref{entries} to give a general condition under which two non-resonant \sq s are equivalent.

\begin{prop} \label{nonres eqvs}
For $l$ arbitrary and $n$ non-resonant with respect to $l$,\/ $\SQ^{\delta-n,l}_\lu$ and\/ $\SQ^{\delta'-n,l}_\lup$ are equivalent \iff\ there are non-zero scalars $\ep_n, \ep_{n+1}, \ldots, \ep_{n+l-1}$ such that for all $(i,j)$ with $0 \le j \le i-2$ and $2 \le i \le l-1$,
\begin{equation} \label{ep}
   (\delta-n-j)_{i-j}\, B_{n+i,n+j}(\gd)\, \ep_{n+i}\ =\ 
   (\delta'-n-j)_{i-j}\, B_{n+i,n+j}(\gdp)\, \ep_{n+j}.
\end{equation}
\end{prop}

\meno {\em Proof.\/}
The two \sq s are equivalent \iff\ the corresponding \r s $\pi^\lu$ and $\pi^\lup$ on $\bigoplus_{i=0}^{l-1} \F_{n+i}$ are equivalent.  Note that both $\pi^\lu$ and $\pi^\lup$ carry the Casimir operator $Q$ to the same block-diagonal matrix with scalars on the diagonal:
\begin{equation*}
   \pi^\lu_{n+i,n+j}(Q) = \pi^\lup_{n+i,n+j}(Q) = \delta_{i,j} (n+i)(n+i-1),
\end{equation*}
where $\delta_{i,j}$ is the Kronecker delta function.  Since $n$ is non-resonant, these scalars are distinct.

Suppose that $\ep$ is an endomorphism of $\bigoplus_{i=0}^{l-1} \F_{n+i}$ intertwining $\pi^\lu$ and $\pi^\lup$.  Since it commutes with the $Q$ action, it must be block-diagonal: write $\ep_n, \ldots, \ep_{n+l-1}$ for its diagonal entries.  Since the diagonal entries $\pi^\lu_{n+i,n+i}$ and $\pi^\lup_{n+i,n+i}$ are both $L_{n+i}$, $\ep_{n+i}$ must intertwine $L_{n+i}$ with itself.  It is elementary that this forces $\ep_{n+i}$ to be a scalar.  By Corollary~\ref{entries} and Theorem~\ref{CMZ}, (\ref{ep}) is the condition for $\ep$ to intertwine $\pi^\lu$ with $\pi^\lup$.  $\Box$

\meno {\em Proofs of Propositions~\ref{nonres3}, \ref{nonres4}, \ref{l2l3l4}, and~\ref{general invnt}, and Theorem~\ref{nonres5}.\/}  These are all essentially corollaries of Proposition~\ref{nonres eqvs}, and we will only discuss Theorem~\ref{nonres5} and Proposition~\ref{general invnt}.  In Theorem~\ref{nonres5}, it is obvious from~(\ref{ep}) that the SVC is necessary for equivalence.  Assuming that it holds, (\ref{ep}) gives one condition on the ratios of the scalars $\ep_{n+i}$ for each of the six quantities~(\ref{6Bijs}) that is non-zero.  Since there are five scalars and hence four independent ratios in length~5, there are in general two additional conditions on $(\gd)$ which must be satisfied in order for~(\ref{ep}) to be soluble.

For example, suppose that none of the quantities~(\ref{6Bijs}) vanishes.  Rewrite~(\ref{ep}):
\begin{equation*}
   \frac{\ep_{n+i}}{\ep_{n+j}}\ =\ 
   \frac{(\delta'-n-j)_{i-j}}{(\delta-n-j)_{i-j}}\,
   \frac{B_{n+i,n+j}(\gdp)}{B_{n+i,n+j}(\gd)}\,.
\end{equation*}
Solving these equations for $(i,j)$ equal to $(2,0)$, $(3,1)$, $(4,2)$, and $(3,0)$ determines the $\ep_{n+i}$ up to a common multiplicative scalar.  The reader may easily check that they are soluble at the two remaining places $(4,1)$ and $(4,0)$ \iff\ Part~(v) of Theorem~\ref{nonres5} holds.

For Proposition~\ref{general invnt}, write $x \cdot B$ and $y \cdot B$ for the numerator and denominator of the displayed equation.  It is easy to check that
\begin{equation*}
   \frac{x \cdot B(\gdp)}{x \cdot B(\gd)}\ =\
   \frac{(\delta - n)_4\, (\delta - n - 1)_2}{(\delta' - n)_4\, (\delta' - n - 1)_2}\,
   \frac{\ep_{n+4}\, \ep_{n+3}}{\ep_{n+1}\, \ep_n},
\end{equation*}
independent of $x$.  The invariance follows, and we leave the reader to check the completeness of $\t I_n$ and $\t J_n$.  $\Box$

\subsection{The resonant case}
Recall from Section~\ref{Dfns} that $n$ is resonant with respect to $l$ if $n \in -\oh\bN$ and $n+l-1 \ge 1$.  In this case some of the eigenvalues of $Q$ on the \cs\ $\{ \F_n, \ldots, \F_{n+l-1} \}$ of $\SQ^{\delta-n,l}_\lu$ are doubled: those which add to~$1$.  

Here there is in general no $\da$-equivalence of the form~(\ref{PQ}), but it is possible to choose an {\em affine quantization\/} $\o\PQ_\lu$ in place of $\PQ_\lu$ which is as close as possible to a projective quantization, in the sense that it preserves the generalized eigenspaces of $Q$.  By an affine quantization we mean a symbol-preserving affine equivalence, where the {\em affine algebra\/} $\db$ is the Borel subalgebra of $\da$ defined by
\begin{equation*}
   \db := \Span_\bC \bigl\{ \px, x\px \bigr\}.
\end{equation*}

Our derivation of the resonant \ec es relies on the explicit form of the \r\ $\o\pi^\lu$ taking the place of $\pi^\lu$ in~(\ref{pilu}).  We now recall from \cite{CS04} the formulas for the matrix entries $\o\pi^\lu_{n+i,n+j}$.  In fact we will give more detail than is needed for the current article, because we take this opportunity to give in Theorem~\ref{res CMZ} a significant simplification of the formulas for the scalars $\o b_{n+i,n+j}(\lu)$ given in Theorem~6.3 of \cite{CS04}.  We begin with the definitions and properties of the two new types of 1-cochains occurring in the entries of $\o\pi^\lu$.  

Write $\partial$ for the coboundary operator: if $T \in \D_\lu$ is 0-cochain and $\omega:\VR \to \D_\lu$ is a 1-cochain, then their coboundaries are, respectively, the 1- and 2-cochains
\begin{equation*}
   \partial T (X) = L_\lu(X)T, \quad
   \partial \omega(X \wedge Y) = L_\lu(X) \omega(Y)
   - L_\lu(Y) \omega(X) - \omega ([X,Y]).
\end{equation*}
A 2-cochain $\theta: \Lambda^2\VR \to \D_\lu$ is said to be {\em $\da$-relative\/} if it is $\da$-covariant and zero on $\da\wedge\VR$.  

Consider the following three properties a 1-cochain $\omega$ may have:
\begin{equation*}
   \mbox{\rm {\em (i)\/} $\omega$ is zero on $\da$,} \qquad
   \mbox{\rm {\em (ii)\/} $\omega$ is $\da$-covariant,} \qquad
   \mbox{\rm {\em (iii)\/} $\partial\omega$ is $\da$-relative.}
\end{equation*}
If $\omega$ has any two of these properties, then it has the third and is $\da$-relative.  However, it is possible to have any one of these properties but neither of the other two.

We will need also the {\em cup product\/} of a $\D_{\nu,\mu}$-valued 1-cochain $\omega_{\nu,\mu}$ with a $\D_{\lambda,\nu}$-valued 1-cochain $\omega_{\lambda,\nu}$:
\begin{equation*}
   \omega_{\nu,\mu} \cup \omega_{\lambda,\nu} (X\wedge Y) :=
   \omega_{\nu,\mu}(X) \circ \omega_{\lambda,\nu}(Y) - 
   \omega_{\nu,\mu}(Y) \circ \omega_{\lambda,\nu}(X).
\end{equation*}

\meno {\bf Definition.}
For $\delta \in 1+\bN$, let $\alpha_\lu: \VR \to \D_\lu$ be
\begin{equation*}
   \alpha_\lu := \frac{2\sqrt{3}}{\delta \gamma^{1/2}}
   \Bigl( {\ts\frac{1}{12}} \delta (\delta-1) (\delta+1-\gamma) \beta_\lu
          - \partial(dx^\delta \px^\delta) \Bigr).
\end{equation*}

Note that although $\beta_\lu$ is not defined at $\delta=1$, it does not appear in the formula there.  For $\delta \in 2+\bN$, $\alpha_\lu$ is manifestly cohomologous to a multiple of $\beta_\lu$ except in the self-adjoint case $\gamma=0$, where it appears to be undefined.  It is a crucial point in \cite{CS04} that in fact the formula for $\alpha_\lu$ has a removable singularity at $\gamma=0$, because this is exactly where it is needed.  At this value, $\beta_\lu$ is a coboundary and the restriction $\alpha_\lu|_\da$ is a non-trivial cocycle.

We define the {\em lacunary \dog\ modules\/} in terms of the lacunary \psdog\ modules of Section~\ref{Lac Results}: for $k \in \bN$ and $k \not= \delta$,
\begin{equation*}
   \D^{k,\lac}_\lu := \D_\lu \cap \Psi^{k,\lac}_\lu.
\end{equation*}

\begin{lemma} \label{alpha}
\begin{enumerate}
\item[(i)]
$\alpha_\lu$ is $\db$-relative, $\D^{\delta-1, \lac}_\lu$-valued, and has symbol
\begin{equation*}
   \alpha_\lu(g\px) \equiv dx^\delta g'' \px^{\delta-1}.
\end{equation*}
\smallbreak \item[(ii)]
For $\delta \not\in 1+\bN$, $C^1\bigl(\VR, \db, \Hom(\F_\lambda, \F_\mu) \bigr) = 0$.
\smallbreak \item[(iii)]
For $\delta=1$, $C^1\bigl(\VR, \db, \Hom(\F_\lambda, \F_\mu) \bigr) = \bC \alpha_\lu$ and $\partial\alpha_\lu = 0$.
\smallbreak \item[(iv)]
For $\delta \in 2+\bN$, $\partial\alpha_\lu$ is $\da$-relative.  Moreover, 
\begin{eqnarray*}
   \bigl\{\omega \in C^1\bigl(\VR, \db, \D^{\delta-1,\lac}_\lu \bigr):
   \mbox{\rm\ $\partial\omega$ is $\da$-relative} \bigr\}
   &=& \bC \alpha_\lu, \mbox{\rm\ and} \\[6pt]
   \bigl\{\omega \in C^1\bigl(\VR, \db,  \Hom(\F_\lambda, \F_\mu) \bigr):
   \mbox{\rm\ $\partial\omega$ is $\da$-relative} \bigr\}
   &=& \Span_\bC\{\alpha_\lu, \beta_\lu\}.
\end{eqnarray*}
\smallbreak \item[(v)]
Under conjugation, $\alpha_\lu^* = (-1)^{\delta-1} \alpha_{1-\mu,1-\lambda}$.
\end{enumerate}
\end{lemma}

Thus $\alpha_\lu$ is $\db$-relative with $\da$-relative coboundary, but is not zero on $\da$.  In the self-adjoint case it is known to be a cocycle \iff\ $1 \le \delta \le 4$.  Let us remark that $\alpha_\lu(g \px)$ is in general not $\PQ_\lu(dx\, g'')$: a formula for the difference is given in Proposition~9.2 of \cite{CS04}.

The second new type of cochain is zero on $\da$ but does not have an $\da$-relative coboundary.  It has two variations:

\meno {\bf Definition.}
For $\nu-\lambda \in 1+\bN$ and $\mu-\nu \in 2+\bN$, let $\Delta_\lu^{\nu,\LST}:\VR\to\D_\lu$ be
\begin{equation*}
   \Delta_\lu^{\nu,\LST} := -2\bigl( \beta_{\nu,\mu} \circ 
   (dx^{\nu-\lambda} \px^{\nu-\lambda}) - \beta_\lu \bigr) / (\nu-\lambda)(\nu+\lambda-1).
\end{equation*}

For $\nu-\lambda \in 2+\bN$ and $\mu-\nu \in 1+\bN$, let $\Delta_\lu^{\nu,\UST}:\VR\to\D_\lu$ be
\begin{equation*}
   \Delta_\lu^{\nu,\UST} := 2\bigl( (dx^{\mu-\nu} \px^{\mu-\nu}) \circ
   \beta_{\lambda,\nu}  - \beta_\lu \bigr) / (\mu-\nu)(\mu+\nu-1).
\end{equation*}

\medbreak
Just as for $\alpha_\lu$, it is a key point in \cite{CS04} that $\Delta_\lu^{\nu,\LST}$ and $\Delta_\lu^{\nu,\UST}$ have removable singularities at $\nu+\lambda=1$ and $\mu+\nu=1$, so we may consider them to be defined at such values.  In that article, $\Delta_\lu^{\nu,\LST}$ is called $\delta_{\mu-\lambda, \nu-\lambda}(\lambda)$ and $\Delta_\lu^{\nu,\UST}$ is called $\delta_{\mu-\nu, \mu-\lambda}(\lambda)$.  The superscripts ``lower'' and ``upper'' refer to the location of these 1-cochains among the matrix entries of $\o\pi^\lu$: as explained in \cite{CS04}, they occur in the {\em lower singular triangle\/} and {\em upper singular triangle,\/} respectively.

\begin{lemma} \label{Delta}
\begin{enumerate}
\item[(i)]
$\Delta_\lu^{\nu,\LST}$ and $\Delta_\lu^{\nu,\UST}$ are $\db$-relative, $\D^{\delta-3}_\lu$-valued, and have symbols
\begin{equation*}
   -\Delta_\lu^{\nu,\LST}(g\px) \equiv \Delta_\lu^{\nu,\UST}(g\px) \equiv \oh dx^\delta g'''' \px^{\delta-3}.
\end{equation*}
\smallbreak \item[(ii)]
$\Delta_\lu^{\nu,\LST}$ is the unique $\D^{\delta-3}_\lu$-valued 1-cochain such that the following 2-cochain is $\da$-relative:
\begin{equation*}
   \partial\Delta_\lu^{\nu,\LST} + \beta_{\nu,\mu} \cup \alpha_{\lambda,\nu} +
   2 (\partial\beta_{\nu,\mu}) \circ (dx^{\nu-\lambda} \px^{\nu-\lambda}) / (\nu-\lambda)(\nu+\lambda-1).
\end{equation*}
\smallbreak \item[(iii)]
$\Delta_\lu^{\nu,\UST}$ is the unique $\D^{\delta-3}_\lu$-valued 1-cochain such that the following 2-cochain is $\da$-relative:
\begin{equation*}
   \partial\Delta_\lu^{\nu,\UST} + \alpha_{\nu,\mu} \cup \beta_{\lambda,\nu} +
   2(dx^{\mu-\nu} \px^{\mu-\nu}) \circ (\partial\beta_{\lambda,\nu}) / (\mu-\nu)(\mu+\nu-1).
\end{equation*}
\smallbreak \item[(iv)]
Conjugation exchanges $\Delta_\lu^{\nu,\LST}$ and $(-1)^{\delta} \Delta_{1-\mu,1-\lambda}^{1-\nu,\UST}$.
\end{enumerate}
\end{lemma}

\meno {\bf Remark.}
Parts~(ii) and~(iii) of this lemma correct an error in Lemma~6.2 of \cite{CS04}, where the summands containing $\partial\beta_{\nu,\mu}$ and $\partial\beta_{\lambda,\nu}$ were omitted.  This error does not affect the other results of \cite{CS04}.

\medbreak
Suppose now that $n$ is resonant.  Then it is not hard to prove that there exists a $\db$-equivalence
\begin{equation*} 
   \o\PQ_\lu: \bigoplus_{i=0}^{l-1} \F_{n+i} \to \SQ^{\delta-n,l}_\lu
\end{equation*}
which preserves symbols and also preserves the generalized eigenspaces of $Q$: for all $0\le i < j \le l-1$ such that $2n+i+j=1$, $\o\PQ_\lu$ maps $\F_{n+i} \oplus \F_{n+j}$ to the $(n+i)(n+i-1)$-generalized eigenspace of $Q$ on $\SQ^{\delta-n,l}_\lu$.

Such equivalences are not unique.  Given one, define a \r\ $\o\pi^\lu$ of\/ $\VR$ on $\bigoplus_{i=0}^{l-1} \F_{n+i}$ by
\begin{equation*}
   \o\pi^\lu(X) := \o\PQ_\lu^{-1} \circ \big(L_\lu(X)\big|_{\SQ^{\delta-n,l}_\lu}\big) \circ \o\PQ_\lu.
\end{equation*}
Regard $\o\pi^\lu$ as an $l \times l$ matrix with entries
\begin{equation*}
   \o\pi^\lu_{n+i,n+j}: \VR \to \Hom(\F_{n+j}, \F_{n+i}).
\end{equation*}
Then just as in Lemma~\ref{lower tri}, $\o\pi^\lu$ is lower triangular with tensor density actions on the diagonal: for $i<j$, $\o\pi^\lu_{n+i,n+j} = 0$, and at $j=i$, $\o\pi^\lu_{n+i,n+i} = L_{n+i}$.  

We are now in position to state the resonant analogs of Corollary~\ref{entries} and Theorem~\ref{CMZ}.  Parts (i)-(iii) of Proposition~\ref{res entries} and Parts~(i) and~(ii) of Theorem~\ref{res CMZ} were proven in \cite{Ga00}, and the remaining parts of the two statements comprise the main results of \cite{CS04}.  In fact, as we mentioned, Parts~(iii) and~(iv) of Theorem~\ref{res CMZ} significantly simplify the formulas obtained in \cite{CS04}, and we will give a proof of this simplification.

\begin{prop} \label{res entries}
$\o\PQ_\lu$ may be chosen so that $\o\pi^\lu$ has the properties listed below.  In the integral resonant case there is a 1-parameter family of such choices, and in the half-integral resonant case there is only one.
\begin{enumerate}
\smallbreak \item[(i)]
For $i=j+1$ and $n+j\not=0$, $\o\pi^\lu_{n+j+1,n+j}=0$.
\smallbreak \item[(ii)]
For $i \ge j+2$, $\o\pi^\lu_{n+i,n+j}$ is given by~(\ref{pi b beta}) in all of the following cases:
\begin{equation*}
   n+j \ge \oh, \quad n+i \le \oh, \quad
   2n+i+j = 0, \quad 2n+i+j = 2.
\end{equation*}
\smallbreak \item[(iii)]
For $i>j$ and $n+i=1-n-j$, there are scalars $a_{1-n-j,n+j}(\lu)$ such that
\begin{equation*}
   \o\pi^\lu_{1-n-j,n+j}\ =\ a_{1-n-j,n+j}(\lu)\, \alpha_{n+j,1-n-j}.
\end{equation*}
\smallbreak \item[(iv)]
For $n+j \le 0$ and $2n+i+j \ge 3$, there are scalars $\o b_{n+i,n+j}(\lu)$ such that
\begin{eqnarray*}
   \o\pi^\lu_{n+i,n+j} &=& \o b_{n+i,n+j}(\lu)\, \beta_{n+j,n+i} \\[2pt]
   && +\, b_{n+i,1-n-j}(\lu)\, a_{1-n-j,n+j}(\lu)\, \Delta_{n+j,n+i}^{1-n-j,\LST}.
\end{eqnarray*}
\smallbreak \item[(v)]
For $n+i \ge 1$ and $2n+i+j \le -1$, there are scalars $\o b_{n+i,n+j}(\lu)$ such that
\begin{eqnarray*}
   \o\pi^\lu_{n+i,n+j} &=& \o b_{n+i,n+j}(\lu)\, \beta_{n+j,n+i} \\[2pt]
   && +\, a_{n+i,1-n-i}(\lu)\, b_{1-n-i,n+j}(\lu)\, \Delta_{n+j,n+i}^{1-n-i,\UST}.
\end{eqnarray*}
\end{enumerate}
\end{prop}

\begin{thm} \label{res CMZ}
\begin{enumerate}
\item[(i)]
The scalar $a_{1,0}(\lu)$ is $-\delta \gamma^{1/2} / 2 \sqrt{3}$.
\smallbreak \item[(ii)]
For $m \in - \oh\bZ^+$, the scalars $a_{1-m,m}(\lu)$ are given by
\begin{equation*}
   a_{1-m,m}\ =\ -\frac{(\delta-m)_{1-2m}\, B_{1-m,m}}{2 \cdot 3^{(1-2m)/2}\, (-2m)!^2}
   \ =\ -\frac{(\delta-m)_{1-2m}\, (\lambda+\mu-1-m)_{1-2m}}{2\, (-2m)!^2}.
\end{equation*}
\smallbreak \item[(iii)]
For $m \in -\oh\bZ^+$ and $2m+r \ge 3$, the scalars $\o b_{m+r,m}(\lu)$ are given by
\begin{align*}
   \o b_{m+r,m}\ =\ & \frac{3^{-r/2}\, (-1)^{2m-1+r}\, (r^2-1)\, (\delta-m)_r}
   {12\, (2m-2+2r)_{r-2}\, (2m-1+r)\, (2m-3+r)!\, (-2m)!} \\[6pt]
   & \times \biggl[ \frac{1}{2} \frac{\partial}{\partial u} \Big|_{u=m}
   \Bigl( B_{u+r,u} - B_{u+r, u+1-2m}\, B_{u+1-2m,u} \Bigr) \\[6pt]
   &\ \ \ \ - \Bigl( \frac{1}{1-2m} + \frac{2r}{r^2-1} -
   \frac{2(2m-1+r)}{(2m-1+r)^2-1} \Bigr)\, B_{m+r,1-m}\, B_{1-m,m} \biggr].
\end{align*}
\smallbreak \item[(iv)]
For $r \ge 3$, the scalars $\o b_{r,0}(\lu)$ are given by
\begin{align*}
   \o b_{r,0}\ =\ \frac{3^{-r/2}\, (-1)^{r-1}\, (\delta)_r}
   {12\, (2r-2)_{r-3}\, (r-3)!}\,
   & \biggl[ \frac{1}{2} \frac{\partial}{\partial u} \Big|_{u=0} 
   \Bigl( B_{u+r,u} - \gamma^{1/2}\, B_{u+r,u+1} \Bigr) \\[6pt]
   &\ - \Bigl( 1 + \frac{2r}{r^2-1} - \frac{2(r-1)}{r^2-2r} \Bigr)\,
   \gamma^{1/2}\, B_{r,1} \biggr].
\end{align*}
\end{enumerate}
All of the scalars $B_{\bullet,\bullet}$ in this theorem are evaluated at $(\gd)$.
\end{thm}

\meno {\it Proof.\/}
Parts~(i) and~(ii) are proven in both \cite{Ga00} and \cite{CS04}.  Parts~(iii) and~(iv) are proven in \cite{CS04}, except that more complicated formulas for $\o b_{m+r,m}$ are given.  After a long and delicate but elementary computation starting from the \cite{CS04} formulas, one arrives at the formulas above, except that instead of the partial derivative given in Part~(iii) one has
\begin{equation*}
   \lim_{\ep \to 0} \ep^{-1}
   \bigl( B_{m+\ep+r,m+\ep} - B_{m+\ep+r, 1-m+\ep}\, B_{1-m+\ep,m+\ep} \bigr)
   \big|_{\lambda, \mu+\ep},
\end{equation*}
and instead of the partial derivative given in Part~(iv) one has
\begin{equation*}
   \lim_{\ep \to 0} \ep^{-1}
   \bigl( B_{\ep+r,\ep} - \gamma^{1/2}\, B_{\ep+r, 1+\ep} \bigr)
   \big|_{\lambda, \mu+\ep},
\end{equation*}
Here evaluation at $(\lambda,\mu+\ep)$ means that instead of using the usual values $3(\lambda+\mu-1)^2$ and $\mu-\lambda$ for $\gamma$ and $\delta$, we use $3(\lambda+\mu-1+\ep)^2$ and $\mu-\lambda+\ep$.

By \cite{CS04}, these limits exist.  Therefore
\begin{equation*}
   B_{m+r,m} = B_{m+r,1-m} B_{1-m,m}
   \mbox{\rm\ \ and\ \ }
   B_{r,0} = \gamma^{1/2} B_{r,1}
\end{equation*}
for all $(\gd)$.  This is not obvious directly from the formula for $B_{\bullet,\bullet}$ and is stated below in Corollary~\ref{B factors}.

Consequently, the two limits displayed above are in fact directional derivatives along the vector $(1,0,1)$ of
\begin{equation*}
   B_{u+r,u} - B_{u+r, u+1-2m}\, B_{u+1-2m,u}
   \mbox{\rm\ \ and\ \ }
   B_{u+r,u} - \gamma^{1/2}\, B_{u+r,u+1},
\end{equation*}
regarded as functions of $(u,\lambda,\mu)$.  For the first function the derivative is evaluated at $(m,\lambda,\mu)$, and for the second it is evaluated at $(0,\lambda,\mu)$.  By Corollary~\ref{B factors}, at these values of $u$ the two functions are identically zero for all values of $\lambda$ and $\mu$, so the result is unchanged if we take the directional derivative along $(1,0,0)$ instead.  The result follows.  $\Box$

\medbreak
Note that the formula for the scalars $\o b$ occurring in Part~(v) of Proposition~\ref{res entries} is not included in Theorem~\ref{res CMZ}.  This is because it may be deduced from Parts~(iii) and~(iv) of Theorem~\ref{res CMZ} together with the following analog of Proposition~\ref{parity}, which is given as Proposition~9.1 in \cite{CS04}.

\begin{prop} \label{res parity}
\begin{enumerate}
\item[(i)]
The parity equations satisfied by the matrices $\pi$ in Parts~(i) and~(ii) of Proposition~\ref{parity} are also satisfied by the matrices $\o\pi$.
\smallbreak \item[(ii)]
The parity equations satisfied by the scalars $b$ in Parts~(iii) and~(iv) of Proposition~\ref{parity} are also satisfied by the scalars $\o b$.  In particular, for $m+r \in 1+\oh\bN$ and $2m+r \le -1$ we have
\begin{equation*}
   \o b_{m+r,m}(\gamma^{1/2}, \delta) =
   (-1)^{r-1} \o b_{1-m,1-m-r}(\gamma^{1/2}, -\delta).
\end{equation*}
\smallbreak \item[(iii)]
The scalars $a$ satisfy the following parity equations:
\begin{equation*}
   a_{1-m,m}(-\gamma^{1/2}, \delta) =
   (-1)^{1-2m} a_{1-m,m}(\gamma^{1/2}, \delta) =
   a_{1-m,m}(\gamma^{1/2}, -\delta).
\end{equation*}
\end{enumerate}
\end{prop}

These results have the following corollary, which is not noticed in \cite{CS04}.  Half of it was already proven in the proof of Theorem~\ref{res CMZ}; the other half follows from Part~(iv) of Proposition~\ref{parity}.  The simplest examples, $B_{3,0}$ and $B_{1,-2}$, may be observed in the formulas displayed below Proposition~\ref{resSVC}.

\begin{cor} \label{B factors}
\begin{enumerate}
\item[(i)]
For $m \in -\oh\bZ^+$ and $2m+r \ge 3$,
\begin{equation*}
   B_{m+r,m} = B_{m+r,1-m} B_{1-m,m}.
\end{equation*}
\smallbreak \item[(ii)]
For $m+r \in 1+\oh\bZ^+$ and $2m+r \le -1$,
\begin{equation*}
   B_{m+r,m} = B_{m+r,1-m-r} B_{1-m-r,m}.
\end{equation*}
\smallbreak \item[(iii)]
For $r \ge 3$, $B_{r,0} = \gamma^{1/2} B_{r,1}$ and $B_{1,1-r} = \gamma^{1/2} B_{0,1-r}$.
\end{enumerate}
\end{cor}

We have now simplified the results of \cite{CS04} sufficiently to give a useful criterion for equivalence of resonant \sq s, analogous to Proposition~\ref{nonres eqvs}.  We first state a lemma we will need; its proof is an easy exercise.

\begin{lemma} \label{b eqvs}
$\Hom_{\db}(\F_\lambda, \F_\mu)$ is\/ $\bC dx^\delta \px^\delta$ if\/ $\delta \in \bN$, and\/ $0$ otherwise.
\end{lemma}

\begin{prop} \label{res eqvs}
For $n$ resonant with respect to $l$,\/ $\SQ^{\delta-n,l}_\lu$ and\/ $\SQ^{\delta'-n,l}_\lup$ are equivalent \iff\ there are non-zero scalars $\ep_n, \ep_{n+1}, \ldots, \ep_{n+l-1}$ and an arbitrary scalar $\zeta$ such that
\begin{enumerate}
\item[(i)]
For $i \ge j+2$, (\ref{ep}) holds in all of the following cases:
\begin{equation*}
   n+j \ge \oh, \qquad n+i \le \oh, \quad
   \mbox{\rm or} \quad 2n+i+j \in \bigl\{0,1,2\bigr\}.
\end{equation*}
\smallbreak \item[(ii)]
If $n$ is integral, then $\delta\, \gamma^{1/2}\, \ep_1\ =\ \delta'\, \gamma'^{\, 1/2}\, \ep_0$.
\smallbreak \item[(iii)]
For $n+j \le -\oh$ and $2n+i+j \ge 3$, or $n+i \ge \frac{3}{2}$ and $2n+i+j \le -1$,
\begin{equation*} 
   \o b_{n+i,n+j}(\gd)\, \ep_{n+i}\ =\ \o b_{n+i,n+j}(\gdp)\, \ep_{n+j}.
\end{equation*}
For $n+j \le -\oh$ and $2n+i+j \ge 3$, this condition is equivalent to
\begin{align*} 
   & \ep_{n+i}\, (\delta-n-j)_{i-j}\,
   \frac{\partial}{\partial u} \Big|_{u=n+j}
   \Bigl( B_{u+r,u} - B_{u+r, u+1-2m}\, B_{u+1-2m,u} \Bigr) \\[6pt]
   =\ & \ep_{n+j}\, (\delta'-n-j)_{i-j}\,
   \frac{\partial}{\partial u} \Big|_{u=n+j}
   \Bigl( B'_{u+r,u} - B'_{u+r, u+1-2m}\, B'_{u+1-2m,u} \Bigr),
\end{align*}
where $B_{\bullet,\bullet}$ is evaluated at $(\gd)$ and $B'_{\bullet,\bullet}$ denotes evaluation at $(\gdp)$.
\smallbreak \item[(iv)]
For $n$ integral, $n+j=0$, and $n+i \ge 3$,
\begin{equation*} 
   \o b_{n+i,0}(\gd)\, \ep_{n+i}\ =\ \o b_{n+i,0}(\gdp)\, \ep_0\,
   +\, b_{n+i,1}(\gdp)\, \zeta.
\end{equation*}
\smallbreak \item[(v)]
For $n$ integral, $n+i=1$, and $n+j \le -2$,
\begin{equation*} 
   \o b_{1,n+j}(\gd)\, \ep_1\, +\, b_{0,n+j}(\gd)\, \zeta
   \ =\ \o b_{1,n+j}(\gdp)\, \ep_{n+j}.
\end{equation*}
\end{enumerate}
\end{prop}

\meno {\it Proof.\/}
We follow the approach taken in proving Proposition~\ref{nonres eqvs}, keeping in mind the formulas of Proposition~\ref{res entries} and Theorem~\ref{res CMZ}.  The two \sq s are equivalent \iff\ the \r s $\o\pi^\lu$ and $\o\pi^\lup$ are equivalent.  Suppose that $\ep$ is an endomorphism of $\bigoplus_{i=0}^{l-1} \F_{n+i}$ intertwining $\o\pi^\lu$ and $\o\pi^\lup$.  Regard it as a block matrix with entries $\ep_{n+i,n+j}: \F_{n+i} \to \F_{n+j}$ as usual.  The fact that $\ep$ preserves the generalized eigenspaces of $Q$ implies that $\ep_{n+i,n+j} = 0$ except possibly on the diagonal $i=j$ and the antidiagonal $2n+i+j=1$.

Let $\ep_n, \ldots, \ep_{n+l-1}$ be the diagonal entries of $\ep$.  Recall that the diagonal entries of the two \r s are $\o\pi^\lu_{m,m} = \o\pi^\lup_{m,m} = L_m$, and that these \r s both arise from affine quantizations, so their restrictions to $\db$ are diagonal.  Therefore by Lemma~\ref{b eqvs} the $\ep_m$ are scalars, and the antidiagonal entries $\ep_{1-m,m}$ are multiples of the {\em Bol operators\/}
\begin{equation*}
   \Bol_{m,1-m}\ :=\ dx^{1-2m}\, \px^{1-2m}: \F_m \to \F_{1-m},
\end{equation*}
which are known to be $\da$-maps.  Observe that $\Bol_{0,1}$ is the de Rham differential $d$, a $\VR$-map.

Consider the intertwining equation $\ep \circ \o\pi^\lu = \o\pi^\lup \circ \ep$.  For $n+j \ge \oh$ or $n+i \le \oh$, its $(n+i, n+j)$ entry is $\o\pi^\lu_{n+i, n+j}\, \ep_{n+i} = \o\pi^\lup_{n+i, n+j}\, \ep_{n+j}$, which yields Part~(i) in these cases.

For $2n+i+j=2$, the $(n+i, n+j)$ entry of the intertwining equation is
\begin{equation*}
   \o\pi^\lu_{2-n-i,n+i}\, \ep_{n+i} = 
   \o\pi^\lup_{2-n-i,n+i}\, \ep_{n+i} + \o\pi^\lup_{2-n-i,1-n-i}\, \ep_{1-n-i,n+i}.
\end{equation*}
Since $\o\pi_{m+1,m} = 0$ for $(m+1,m) \not= (1,0)$, we obtain Part~(i) in this case also.  The proof of Part~(i) for $2n+i+j=0$ is similar.

On the antidiagonal $2n+i+j=1$, recall that the $(1-m, m)$-entry $\o\pi^\lu_{1-m,m}$ is $a_{1-m,m} \alpha_{1-m,m}$.  Use the corresponding entry of the intertwining equation to obtain the proportionality
\begin{equation*}
   \bigl( \ep_{1-m}\, a_{1-m,m}(\lu) - \ep_m\, a_{1-m,m}(\lup) \bigr) \alpha_{1-m,m}
   =  \partial \ep_{1-m,m} \propto \partial \Bol_{m,1-m}.
\end{equation*}
It is proven in \cite{CS04} that $\partial \Bol_{m,1-m} = \frac{1}{12} (2-2m)_2\, \beta_{m,1-m}$, a non-zero multiple of $\beta_{1-m,m}$ unless $m=0$.  (Indeed, this is why $\alpha_\lu$ is defined at $\lambda + \mu = 1$.)  Since $\alpha_{1-m,m}$ and $\beta_{1-m,m}$ are linearly independent, both sides of the above displayed equation must be zero.  This completes the proof of Part~(i) and also proves Part~(ii).

The argument in the last paragraph also gives $\ep_{1-m,m} = 0$ for $m \not= 0$, and, if $n$ is integral, $\ep_{1,0} = \zeta d$ for some scalar $\zeta$.  Since $\ep_{1,0}$ does not enter the intertwining equation unless $n+j=1$ or $n+i=0$, the first sentence of Part~(iii) now follows easily: the terms involving the cochains $\Delta_{n+j,n+i}$ cancel due to the conditions imposed by Part~(i).  For the second sentence of Part~(iii), use Part~(iii) of Theorem~\ref{res CMZ}: the terms not involving $\partial_u$ cancel due to Part~(i) of the current proposition.  

Parts~(iv) and~(v) are similar so we only discuss Part~(iv).  For $r \ge 3$, the $(r,0)$-entry of the intertwining equation is
\begin{equation*}
   \ep_r\, \o\pi^\lu_{r,0}\, =\, \ep_0\, \o\pi^\lup_{r,0}
   +\, \zeta\, \o\pi^\lup_{r,1}\, d.
\end{equation*}
Again, the terms involving the cochains $\Delta_{n+j,n+i}$ cancel due to the conditions imposed by Parts~(i) and~(ii).  To finish the proof, use Lemma~\ref{beta props} together with the fact that $d$ is an $\da$-map to see that $\beta_{r,1} \circ d = \beta_{r,0}$.  $\Box$

\meno {\it Proofs of Propositions~\ref{resSVC}, \ref{res l=4 n=0}, and~\ref{resSD5}, and Theorem~\ref{resSD4}.\/}  Just as the non-resonant equivalence results are all corollaries of Proposition~\ref{nonres eqvs}, the resonant results are all corollaries of Proposition~\ref{res eqvs}.  We shall only give the details in the most important cases.

In Proposition~\ref{resSD5}, all of the entries of the intertwining equation involved fall under Part~(i) of Proposition~\ref{res eqvs}.  Therefore the proof goes exactly as for Theorem~\ref{nonres5}.

The proof of Theorem~\ref{resSD4} is similar: only Parts~(i) and~(ii) of Proposition~\ref{res eqvs} are involved.  In the event that none of the four functions in the theorem vanishes, solving~(\ref{ep}) for $(n+i, n+j)$ equal to $(1, -1)$, $(2, -1)$, and $(1, 0)$ determines the scalars $\ep_{-1}$, $\ep_0$, $\ep_1$, and $\ep_2$ up to a common multiplicative scalar.  Then solving it at the last entry $(2, 0)$ gives the invariant $R$.  The scalar $\zeta$ has no effect.

The $n = -2$ case of Proposition~\ref{res l=4 n=0} follows from the $n=0$ case by duality.  The $n=0$ case involves Part~(iv) of Proposition~\ref{res eqvs}, but if $b_{3,1} \not= 0$ we can pick $\zeta$ to satisfy the relevant equation, imposing no new condition.  When $b_{3,1} = 0$, $\zeta$ has no effect and we get the additional condition
\begin{equation*}
   \o b_{3,0}(\gd)\, \ep_3 = \o b_{3,0}(\gdp)\, \ep_0.
\end{equation*}
Using Theorem~\ref{res CMZ}, one obtains $\o b_{3,0}\, =\, -\frac{1}{2\sqrt{3}}\, (\delta)_3\, \gamma^{1/2}\, \bigl( \frac{5}{6} B_{3,1} + \delta - 3 \bigr)$.  Proposition~\ref{res l=4 n=0} now follows easily.  $\Box$

\subsection{Lacunary cases}
First note that the $\VR$-module $\Psi^{k,\lac}_\lu$ exists because the entries $\pi_{m+1,m}$ and $\o\pi_{m+1,m}$ on the first subdiagonal are all zero except for $\o\pi_{1,0}$.  Each lacunary module still has either a projective quantization $\PQ_\lu$ or an affine quantization $\o\PQ_\lu$, depending on whether or not it is resonant.  To obtain the associated \r\ $\pi$ or $\o\pi$, simply delete from the non-lacunary $\pi$ or $\o\pi$ the rows and columns passing through the diagonal entries $\pi_{m,m}$ or $\o\pi_{m,m}$ corresponding to the excised composition series modules $\F_m$.

We will only give the proofs of Theorems~\ref{lac024} and~\ref{lac0235}; the other lacunary results are proven similarly.

\meno{\it Proof of Theorem~\ref{lac024}.\/}
In the non-resonant case we obtain the \r s
\begin{equation*}
   \pi^\lu = \begin{pmatrix} L_n & & \\
   \pi^\lu_{n+2,n} & L_{n+2} & \\
   \pi^\lu_{n+4,n} & \pi^\lu_{n+4,n+2} & L_{n+4} \end{pmatrix}, \quad
   \pi^\lup = \begin{pmatrix} L_n & & \\
   \pi^\lup_{n+2,n} & L_{n+2} & \\
   \pi^\lup_{n+4,n} & \pi^\lup_{n+4,n+2} & L_{n+4} \end{pmatrix},
\end{equation*}
both acting on $\F_n \oplus \F_{n+2} \oplus \F_{n+4}$.  The subdiagonal entries of these matrices are still given by the formulas of Corollary~\ref{entries} and Theorem~\ref{CMZ}.  The reasoning of Proposition~\ref{nonres eqvs} shows that any equivalence $\ep$ between them must be block-diagonal with scalar diagonal entries $\ep_n$, $\ep_{n+2}$, and $\ep_{n+4}$.  If none of the subdiagonal entries is zero, then solving~(\ref{ep}) at the entries $(n+2,n)$ and $(n+4,n+2)$ determines these three scalars up to a common multiplicative scalar, and solving it at $(n+4,n)$ gives the invariant $I_n$.  This proves the theorem in the non-resonant cases.  The only resonant cases involve only Part~(i) of Theorem~\ref{res eqvs}, so the non-resonant proof is valid also for them.  $\Box$

\meno{\it Proof of Theorem~\ref{lac0235}.\/}
Here the \r\ $\pi^\lu$ is
\begin{equation*}
   \begin{pmatrix} L_n & & & \\
   \pi^\lu_{n+2,n} & L_{n+2} & & \\
   \pi^\lu_{n+3,n} & 0 & L_{n+3} & \\
   \pi^\lu_{n+5,n} & \pi^\lu_{n+5,n+2} & \pi^\lu_{n+5,n+3} & L_{n+5}
   \end{pmatrix}.
\end{equation*}
If none of the subdiagonal entries is zero, then solving~(\ref{ep}) at the entries $(n+2,n)$, $(n+3,n)$, and $(n+5,n+3)$ determines the scalars $\ep_n$, $\ep_{n+2}$, $\ep_{n+3}$, and $\ep_{n+5}$ up to a common multiplicative scalar, and solving it at $(n+5,n+2)$ gives the invariant $M_n$.

At this point a new phenomenon occurs.  One would expect~(\ref{ep}) to impose a new condition at $(n+5,n)$, but since $n$ is not allowed to be $-4$ or $0$, the following proposition shows that it does not: (\ref{ep})~is automatically satisfied at $(n+5, n)$ if it is satisfied at the other entries.  Therefore the proof is complete.  $\Box$

\begin{prop} \label{noncocycles}
Fix complex scalars $n_1, \ldots, n_l$ such that $n_i - n_{i-1} \in \bZ^+$ for $1 \le i \le l$.  Suppose that $\sigma$ and $\sigma'$ are two \r s of\/ $\VR$ on $\bigoplus_{i=1}^l \F_{n_i}$.  Regard them as block matrices with entries $\sigma_{i,j}: \VR \to \Hom(\F_{n_j}, \F_{n_i})$.  Assume that $\sigma$ has the following form:
\begin{enumerate}
\item[(i)]
$\sigma$ is lower triangular: $\sigma_{i,j} = 0$ for $i < j$.
\smallbreak \item[(ii)]
The diagonal entries are the Lie actions: $\sigma_{i,i} = L_{n_i}$.
\smallbreak \item[(iii)]
For $n_{j+1} - n_j = 1$ and $n_j \not= 0$, $\sigma_{j+1,j} = 0$.
\smallbreak \item[(iv)]
For $n_i - n_j \ge 2$ and $n_i + n_j \not= 1$, there is a scalar $c_{i,j}$ such that
\begin{equation*}
   \sigma_{i,j}\, =\, c_{i,j}\, \beta_{n_j,n_i}.
\end{equation*}
\smallbreak \item[(v)]
For $i > j$ and $n_i + n_j = 1$, there is a scalar $c_{i,j}$ such that
\begin{equation*}
   \sigma_{i,j}\, =\, c_{i,j}\, \alpha_{n_j,n_i}.
\end{equation*}
\end{enumerate}
Assume that $\sigma'$ has the same form with scalars $c'_{i,j}$ replacing the $c_{i,j}$.

Suppose that $\ep$ is an endomorphism of\/ $\bigoplus_{i=1}^l \F_{n_i}$ whose block matrix is diagonal, with non-zero scalars $\ep_1, \ldots, \ep_l$ on the diagonal satisfying
\begin{equation*}
   c_{i,j}\, \ep_i\, =\, c'_{i,j}\, \ep_j
\end{equation*}
for all $(i,j)$ such that either $n_i - n_j \in \{2, 3, 4\}$ or $(n_i, n_j)$ is one of the following:
\begin{equation} \label{cocycles}
   (1,0), \quad (1, -4), \quad (5,0), \quad \oh(7 \pm \sqrt{19}, -5 \pm \sqrt{19}).
\end{equation}
Then $c_{i,j} \ep_i = c'_{i,j} \ep_j$ for all $(i,j)$, and $\ep$ is an intertwining map from $\sigma$ to $\sigma'$.
\end{prop}

\meno {\it Proof.\/}
Replacing $\sigma$ with $\ep \circ \sigma \circ \ep^{-1}$, we reduce to proving that the entries $\sigma_{i,j}$ with $2 \le n_i - n_j \le 4$ or $(n_i, n_j)$ as in~(\ref{cocycles}) determine $\sigma$.  Recall the following elementary cohomological result: the condition that $\sigma$ be a \r\ is equivalent to the {\it cup equation:\/}
\begin{equation*}
   0\, =\, \partial \sigma_{i,j} +
   \sum_{i>k>j} \sigma_{i,k} \cup \sigma_{k,j} \quad
   \mbox{\rm for all\ } i > j.
\end{equation*}
Thus $c_{i,j}$ is determined by the entries on the higher subdiagonals whenever $\partial \beta_{n_j,n_i}$ (or $\partial \alpha_{n_j,n_i}$ if $n_i + n_j = 1$) is non-zero.  The result now follows from the well-known theorem of Feigin and Fuchs on the cohomology of $\VR$ \cite{FF80}, of which a small part tells us that $\partial \beta_{\lambda,\mu} = 0$ \iff\ $\mu - \lambda \in \{2,3,4\}$ or $(\lambda,\mu)$ is one of $(-4,1)$, $(0,5)$, or $\oh(-5 \pm \sqrt{19}, 7 \pm \sqrt{19})$, and $\partial \alpha_{\lambda, 1-\lambda} = 0$ \iff\ $\lambda \in \{0, -\oh, -1, -\frac{3}{2}\}$.  $\Box$

\section{Remarks} \label{Remarks}

We conclude with a discussion of the state of the equivalence question.  The main open problem is to answer it in lengths $l > 5$.  In the resonant case one could try to do this with Proposition~\ref{res eqvs} and the formulas of Theorem~\ref{res CMZ}, but we have not taken any steps in this direction.  In the non-resonant case the necessary tools are Proposition~\ref{nonres eqvs} and the formula of Theorem~\ref{CMZ}.  Consider the case of length~6.  The first point is that by Proposition~\ref{noncocycles}, the $(n+5,n)$ entry of~(\ref{ep}) does not impose any new condition on $\ep$, so we have the following proposition.

\begin{prop}
For $n$ non-resonant, the \sq s $\SQ^{\delta-n,6}_\lu$ and\/ $\SQ^{\delta'-n,6}_\lup$ are equivalent \iff\ the two pairs of length~5 \sq s $\bigl( \SQ^{\delta-n,5}_\lu,\, \SQ^{\delta'-n,5}_\lup \bigr)$ and\/ $\bigl( \SQ^{\delta-n-1,5}_\lu,\, \SQ^{\delta'-n-1,5}_\lup \bigr)$ are both equivalent.
\end{prop}

Thus in the generic case where both pairs of length~5 \sq s fall under Part~(v) of Theorem~\ref{nonres5}, the four functions $I_n$, $J_n$, $I_{n+1}$, and $J_{n+1}$ are complete invariants for the length~6 \ec es $\EC^6_n(\gd)$.  Similarly, in the non-resonant length~7 case the entries at $(n+5,n)$ and $(n+6,n+1)$ do not impose new conditions, and the entry at $(n+6,n)$ can do so only when the Feigin-Fuchs cocycles are involved, at $2n = -5 \pm \sqrt{19}$.  It would be very interesting if these special values of $n$ arise in the answer to the length~7 equivalence question, but we expect that they do not: we conjecture that in the non-resonant case, the only equivalences in lengths $l \ge 6$ are given by conjugation and the de Rham differential, as explained in Lemmas~\ref{conj} and~\ref{Bol}, respectively.  Despite substantial effort and use of computers we have been unable to prove this conjecture in length~6 or~7, but it is relatively easy to prove in length~8 and hence in all higher lengths.  We simply state the result here; we plan to give the proof in a future article treating also the length~6 and~7 cases.

\begin{prop}
For $n$ non-resonant and $l \ge 8$, $\SQ^{\delta-n,l}_\lu$ and $\SQ^{\delta'-n,l}_\lup$ are equivalent \iff\ their parameters $(\gd)$ and $(\gdp)$ are either equal, in which case the modules are either conjugate or themselves equal, or make up the pair $\bigl( 3(\nu+1)^2, \nu \bigr)$, $\bigl( 3\nu^2, \nu+1 \bigr)$ for some $\nu$, the case of Lemma~\ref{Bol}.
\end{prop}

We remark that analysis of the \ec es of the lacunary modules $\Psi^{\delta-n,\lac}_\lu / \Psi^{\delta-n-5,\lac}_\lu$ may be seen as an interesting ``warm-up'' problem for the equivalence question in length~7, as they are \sq s of $\SQ^{\delta-n,7}_\lu$.  These modules have \cs\ $\{ \F_n, \F_{n+2}, \F_{n+3}, \F_{n+4}, \F_{n+6} \}$, and by Proposition~\ref{noncocycles} there are three rational invariants for their \ec es: $I_n$, $I_{n+2}$, and
\begin{equation*}
   B_{n+6,n+3} B_{n+3,n} \big/ B_{n+6,n+4} B_{n+4,n+2} B_{n+2,n},
\end{equation*}
whose level curves form a pencil of cubics.  Coupled with the SVC these invariants are complete except possibly at $2n = -5 \pm \sqrt{19}$, where one may need also the invariant
\begin{equation*}
   B_{n+6,n} / B_{n+6,n+3} B_{n+3,n},
\end{equation*}
whose level curves form another pencil of cubics.

\def\eightit{\it} 
\def\bib{\bf}
\bibliographystyle{amsalpha}

\end{document}